\documentclass[12pt,reqno]{amsart}

\usepackage[pdftex]{graphicx} 
\usepackage[english]{babel} 
\usepackage[pdftex,linkcolor=black,pdfborder={0 0 0}]{hyperref} 
\usepackage{calc} 
\usepackage{enumitem} 
\usepackage{comment}

\usepackage{import}

\usepackage{adjustbox}

\usepackage{amsmath}
\usepackage{amssymb}
\usepackage{amsthm}
\usepackage[all]{xy}
\usepackage[pdftex]{graphicx}
\usepackage{color}
\usepackage{cite}
\usepackage{url}
\usepackage{indent first}
\usepackage[labelfont=bf,labelsep=period,justification=raggedright]{caption}
\usepackage[english]{babel}
\usepackage[utf8]{inputenc}
\usepackage[colorinlistoftodos]{todonotes}
\usepackage{tkz-fct}
\usepackage{tikz}
\usetikzlibrary{calc}
\usepackage{comment}
\usepackage{fullpage}

\usepackage{mathtools} 
\mathtoolsset{showonlyrefs} 

\usepackage{multicol}
\PassOptionsToPackage{dvipsnames,svgnames}{xcolor}
\usepackage{textcomp}
\usepackage{sistyle}
\SIthousandsep{,}

\usepackage{fullpage}

\usepackage{lipsum} 
\usepackage{import}

\usepackage{tikz}
\usetikzlibrary{decorations.markings}

\newcommand{\Q}{\mathbb{Q}}

\newcommand{\K}{\mathbb{K}}
\newcommand{\Z}{\mathbb{Z}}

\newcommand{\PP}{\mathfrak{p}}

\renewcommand{\O}{\mathcal{O}}

\theoremstyle{plain}
\newtheorem{theorem}{Theorem}[section]

\newtheorem{lemma}[theorem]{Lemma}

\theoremstyle{definition}

\newtheorem*{remark}{Remark}

\newcommand{\erf}[1]{\text{erf}\left(#1\right)}
\newcommand{\erfc}[1]{\text{erfc}\left(#1\right)}

\renewcommand{\Re}{\operatorname{Re}}

\newcommand{\red}[1]{{\color{red} #1}}

\begin{document}

\title[The number of integral ideals]{On the number of integral ideals in a number field}

\author[E.~S.~Lee]{Ethan~Simpson~Lee}
\address{University of Bristol, School of Mathematics, Fry Building, Woodland Road, Bristol, BS8 1UG} 
\email{ethan.lee@bristol.ac.uk}
\urladdr{\url{https://sites.google.com/view/ethansleemath/home}}

\maketitle

\begin{abstract}
We update Sunley's explicit estimate for the ideal-counting function, which is the number of integral ideals of bounded norm in a number field.
\end{abstract}

\section{Introduction}

Suppose that a number field $\K$ has degree $n_{\K}$, ring of integers $\mathcal{O}_{\K}$, and the absolute value of its discriminant is $| \Delta_{\K} |$. The Dedekind zeta-function associated to $\K$, denoted $\zeta_{\K}(s)$, is regular throughout $\mathbb{C}$ aside from one pole at $s=1$ which is simple and has residue $\kappa_{\K}$. Throughout this paper, we will use big-$O$ and $\ll$ notation, in which the implied constants depend on $\K$.

The number $I_{\K}(x)$ of ideals $\mathfrak{a}\subset \O_{\K}$ with norm $N(\mathfrak{a})\leq x$ is called the ideal-counting function. Estimates for $I_{\K}(x)$ are naturally useful, because it is the number fields generalisation of the \textit{floor} (or \textit{integer counting}) function, which is commonly denoted $[x]$ or $\lfloor x \rfloor$. In this paper, we investigate the error term in the well-known estimate $I_{\K}(x) \sim \kappa_\K x$  as $x\to\infty$. To this end, Weber \cite{Weber} showed
\begin{equation}\label{eqn:Landau_similar}
    I_{\K}(x) = \kappa_\K x + O\!\left(x^{1-\frac{1}{n_\K}}\right).
\end{equation}
Soon after, Landau \cite[Satz 210]{LandauEinfuehrung} improved this to
\begin{equation}\label{eqn:Landau}
    I_{\K}(x) = \kappa_\K x + O\!\left(x^{1 - \frac{2}{n_{\K} + 1}}\right).
\end{equation}
Shapiro used an estimate of the form \eqref{eqn:Landau_similar} in \cite{Shapiro}, to generalise Erd\H{o}s and Selberg's elementary proof of the prime number theorem \cite{Erdos, Selberg}. Moreover, \eqref{eqn:Landau_similar} and \eqref{eqn:Landau} are important in the methods demonstrated in \cite{Rosen} to establish Mertens' theorems for number fields.

There has been a lot of work done to improve the shape of the error terms in \eqref{eqn:Landau_similar} and \eqref{eqn:Landau}. Using the Hardy--Littlewood circle method, Huxley and Watt improved \eqref{eqn:Landau} for quadratic fields in \cite{HuxleyWatt}, M\"{u}ller tackled the cubic case in \cite{Muller}, and Bordell\`{e}s has improved \eqref{eqn:Landau} whenever $n_{\K}\geq 4$ in \cite{Bordelles}. Bordell\`{e}s' work generalised the work of M\"{u}ller \cite{Muller}, and built upon Nowak \cite{Nowak}.

Another direction of research seeks explicit statements for \eqref{eqn:Landau_similar} and \eqref{eqn:Landau}. Here, explicit means that the implied constant in the error term is fully described using explicit constants, depending only on the invariants $n_{\K}$ and $\Delta_{\K}$ of the number field $\K$. To this end, Theorem \ref{thm:Sunley} (below) is an explicit version of \eqref{eqn:Landau} which was established in Sunley's thesis \cite{SunleyThesis}; it is also presented in \cite[Thm.~2]{SunleyBAMS} and \cite[Thm.~1.1]{SunleyTAMS} without proof. One can further note that Debaene established an explicit version of \eqref{eqn:Landau_similar} in \cite[Cor.~2]{Debaene}.

\begin{theorem}[Sunley]\label{thm:Sunley}
For $x > 0$ and $n_{\K} \geq 2$, we have
\begin{equation}\label{eqn:Sunley}
    | I_{\K}(x) - \kappa_{\K} x | \leq \Lambda_S(n_{\K}) |\Delta_{\K}|^{\frac{1}{n_{\K}+1}}  (\log | \Delta_{\K} |)^{n_{\K}} x^{1 - \frac{2}{n_{\K}+1}},
\end{equation}
in which $\Lambda_S(n_{\K}) = e^{28.2n_{\K}+5} (n_{\K}+1)^{\frac{5(n_{\K}+1)}{2}}$.
\end{theorem}

Note that it is not prohibitive to restrict our attentions to $n_{\K} \geq 2$, because $I_{\K}(x) = \lfloor x\rfloor$ if $n_{\K} = 1$, and this is well understood. Theorem \ref{thm:Sunley_improved} refines Theorem \ref{thm:Sunley} and is the main result of this paper.

\begin{theorem}\label{thm:Sunley_improved}
For $x > 0$ and $n_{\K} \geq 2$, we have
\begin{equation}\label{eqn:Sunley}
    | I_{\K}(x) - \kappa_{\K} x | \leq \Lambda_{\K}(n_{\K}) |\Delta_{\K}|^{\frac{1}{n_{\K}+1}}  (\log{|\Delta_{\K}|})^{n_{\K}-1} x^{1 - \frac{2}{n_{\K}+1}},
\end{equation}
in which
\begin{equation*}
    \Lambda_{\K}(n_{\K}) = \frac{0.54 (3 n_{\K}-1) \lambda_{\K}(n_{\K})}{(n_{\K}-1)^2 (\log{\mathfrak{m}(n_{\K})})^{n_{\K}-1}} n_{\K}^{3/2} n_{\K}!,
\end{equation*}
$\mathfrak{m}(n_{\K}) = \left(\frac{\pi}{4}\right)^{n_{\K}} n_{\K}^{2 n_{\K}} n_{\K}!^{-2}$, and
\begin{equation*}
    \lambda_{\K}(n_{\K}) =
    \begin{cases}
        (n_{\K}+1)^{\frac{1}{2} - \frac{1}{2 n_{\K}}} \left(\frac{5}{8} + \frac{\pi}{2} - \frac{1}{n_{\K}} + \frac{3}{8 n_{\K}^2}\right)^{\frac{1}{2}} e^{n_{\K}\left(2.27 + \frac{4 n_{\K}}{n_{\K}-1} + \frac{0.01}{n_{\K}^2} + \frac{1}{500 n_{\K}^6}\right)} &\text{if }n_{\K} \leq 13, \\
        (n_{\K}+1)^{n_{\K}-\frac{1}{2}-\frac{1}{2\,n_{\K}}} \left(\frac{5}{8} + \frac{\pi}{2} + \frac{1}{n_{\K}} + \frac{3}{8 n_{\K}^2}\right)^{\frac{1}{2}} e^{4.13 n_{\K} + \frac{0.02}{n_{\K}}} &\text{if }n_{\K} > 13.
    \end{cases}
\end{equation*}
\end{theorem}

An automatic consequence of Theorem \ref{thm:Sunley_improved} is an improvement to the explicit Mertens' theorems for number fields established in \cite{GarciaLeeRamanujan}. Using this, one can also refine \cite[Thm.~2]{GarciaLeeSuhYu}, which is an effective, analytic formula for the number of distinct, irreducible factors of a polynomial. That result shows that there is a finite list (with an effective upper bound) of primes which certifies $f \in \mathbb{Z}[x]$ has exactly $k$ distinct, irreducible factors, and the smaller $\Lambda_{\K}(n_{\K})$ is, the shorter the list of certifying primes will need to be.

We compare $\Lambda_S(n_{\K})$ against $\Lambda_{\K}(n_{\K})$ in Table \ref{tab:SunleyImprovedComparison}. To prove Theorem  \ref{thm:Sunley_improved}, we follow Sunley's proof of Theorem \ref{thm:Sunley} (which is an explicit version of Landau's proof of \eqref{eqn:Landau}). Our improvements come from several avenues. Notably, we optimise over new parameters that we introduce into the method, we update the key ingredients in the method, and we implement modern knowledge about the invariants of the Dedekind zeta-function $\zeta_{\K}(s)$. These key ingredients are explicit upper bounds for $I_{\K}(x)$, $|K(w)|$, and $|K^{(n_{\K})}(w)|$, where $K(w)$ and $K^{(n_{\K})}(w)$ are defined in \eqref{eqn:Kw_defn} and \eqref{eqn:Knw_defn}. The updated ingredients are presented in Theorems \ref{thm:Magdeline} and \ref{thm:Gertrude}.

\begin{table}[]
    \centering
    \begin{tabular}{ccc}
        $n_{\K}$ & $\Lambda_S(n_{\K}) $ & $\Lambda_{\K}(n_{\K})$ \\
        \hline
        2 & $1.75425\cdot 10^{30}$ & $2.49133\cdot 10^{10}$ \\
        3 & $8.57799\cdot 10^{44}$ & $8.45088\cdot 10^{11}$ \\
        4 & $7.88887\cdot 10^{59}$ & $9.84482\cdot 10^{13}$ \\
        5 & $1.20023\cdot 10^{75}$ & $1.41763\cdot 10^{16}$ \\
        10 & $1.90904\cdot 10^{153}$ & $9.65555\cdot 10^{26}$ \\
        15 & $1.10367\cdot 10^{234}$ & $5.27930\cdot 10^{38}$
    \end{tabular}
    \caption{Comparison between values of $\Lambda_S(n_{\K}) $ in Theorem \ref{thm:Sunley} and $\Lambda_{\K}(n_{\K}) $ in Theorem \ref{thm:Sunley_improved} for several choices of $n_{\K}$.}
    \label{tab:SunleyImprovedComparison}
\end{table}

\subsection*{Future work}

In \cite[Thm.~3.3.5]{SunleyThesis}, Sunley proves an explicit P\'{o}lya--Vinogradov theorem for number fields, using analogous methods and ingredients as in her proof of Theorem \ref{thm:Sunley}. In the future, an interested reader could apply analogous ideas or concepts to those laid out in this paper, to improve her explicit P\'{o}lya--Vinogradov theorem for number fields.

\subsection*{Structure}

The goal of this paper is to prove Theorem \ref{thm:Sunley_improved}. In Section \ref{sec:prelim_results}, we introduce several preliminary observations that will be used throughout. In Section \ref{sec:MainTheoremProof}, we prove Theorem \ref{thm:Sunley_improved} using Theorem \ref{thm:Magdeline} and Theorem \ref{thm:Gertrude} in an explicit version of Landau's method to establish \eqref{eqn:Landau}. The key ingredients Theorem \ref{thm:Magdeline} and Theorem \ref{thm:Gertrude} are proved in Sections \ref{sec:IKxUpperBound} and \ref{sec:KwKnwUpperBounds} separately, because they are quite technical.

\subsection*{Acknowledgements}

I would like to thank Richard Brent, Michaela Cully-Hugill, Stephan Garcia, Olivier Ramar\'{e}, Timothy Trudgian, the referee, and all my other colleagues for helpful discussions and comments given throughout the production of this paper.

\subsection*{Quick notes on this arXiv update}

Since the first version of this paper was released, the author was able to refine the main result of the last version into Theorem \ref{thm:Sunley_improved}. This refined result is a proper reflection of the result in the published version of this paper (see \cite{LeeICF}) and the author's PhD thesis \cite{LeeThesis}. Furthermore, there is an argument (in Section \ref{sec:MainTheoremProof}) that should have been clarified in \cite{LeeICF}. For completeness, we include this minor (and inconsequential) update in this article and highlight any differences to the published version \cite{LeeICF} in \red{red}.

\section{Preliminaries}\label{sec:prelim_results}

Throughout, we say that the number field $\K$ has degree $n_{\K} = r_1 + 2r_2$, in which $r_1$ is the number of real places and $r_2$ is the number of complex places of $\K$. Further, suppose $r = r_1 + r_2 - 1$, $\Delta_\K$ is the discriminant of $\K$, $R_{\K}$ is the regulator of $\K$, and $h_{\K}$ is the class number of $\K$. 

\subsection{The Dedekind zeta-function}

Landau establishes most of the knowledge we will state here in \cite{LandauEinfuehrung}. The Dedekind zeta-function is denoted and defined for $\Re{s} > 1$ by
\begin{equation*}
    \zeta_{\K}(s) = \sum_{\mathfrak{a}} N(\mathfrak{a})^{-s} = \prod_{\PP} \left(1 - N(\PP)^{-s}\right)^{-1},
\end{equation*}
which converges absolutely. Here, $\mathfrak{a}\subset \O_{\K}$ are integral ideals of $\K$ and $\PP\subset \O_{\K}$ are \textit{prime} ideals of $\K$. Now, $\zeta_{\K}(s)$ is regular for all $s\in\mathbb{C}$, aside from one simple pole at $s=1$ whose residue is
\begin{equation*}
    \kappa_{\K} = \frac{2^{r_1 + r_2}\pi^{r_2}h_{\K}R_{\K}}{w_{\K}|\Delta_{\K}|^{\frac{1}{2}}};
\end{equation*}
this explicit relationship is called the class analytic formula. For $\Re{s} > 1$, one could alternately re-write the definition of $\zeta_{\K}(s)$ as 
\begin{equation}\label{eqn:DZF_new_def}
    \zeta_{\K}(s) = \sum_{m=1}^\infty \frac{\tau_m}{m^{s}},
\end{equation}
in which $\tau_m$ denotes the number of integral ideals $\mathfrak{a}\subset\O_{\K}$ such that $N(\mathfrak{a}) = m$. The functional equation is $\xi_{\K}(s) = \xi_{\K}(1-s)$, where
\begin{equation}\label{eqn:functional_equation}
    \xi_{\K}(s) = {A_{\K}}^s \Gamma\left(\frac{s}{2}\right)^{r_1} \Gamma\left(s\right)^{r_2} \zeta_{\K}(s)
    \quad\text{such that}\quad
    A_{\K} = 2^{-r_2} \pi^{-\frac{n_{\K}}{2}} |\Delta_{\K}|^{\frac{1}{2}}.
\end{equation}
Using this functional equation, one can deduce $\zeta_{\K}(1-s) = f(s) \zeta_{\K}(s)$ such that
\begin{equation}\label{eqn:LandauSatz156}
    f(s) 
    = \left(\frac{2^{1-s}}{\pi^s}\right)^{n_{\K}} |\Delta_{\K}|^{s-\frac{1}{2}}\left(\cos{\frac{\pi s}{2}}\right)^{r_1 + r_2} \left(\sin{\frac{\pi s}{2}}\right)^{r_2} \Gamma(s)^{n_{\K}}.
\end{equation}

At $s=0$, $\zeta_{\K}(s) = 0$ as long as $r > 0$ and this zero at $s=0$ has order $r$. Further, $\zeta_{\K}(s) = 0$ whenever $s$ is a negative, \textit{even} integer (these zeros have order $r_1+r_2$) or $s$ is a negative, \textit{odd} integer (these zeros only occur when $r_2>0$ and they have order $r_2$). Alongside the zero at $s=0$ (whenever $r>0$), these zeros are called \textit{trivial}. The \textit{non-trivial} zeros of $\zeta_{\K}(s)$ satisfy $0 < \Re{s} < 1$, and we note that there might exist a single, simple, real zero $0 < \beta_0 < 1$, which is called an \textit{exceptional} zero. Explicit bounds for $\beta_0$ may be found in \cite{Ahn, Kadiri, Lee}.

We will require some results moving forward. First, Louboutin has shown in \cite{Louboutin00} that
\begin{equation}\label{eq:LouboutinKappa}
    \kappa_{\K} \leq \left(\frac{e\log{|\Delta_{\K}|}}{2(n_{\K} - 1)}\right)^{n_{\K} - 1}
    \quad \text{for} \quad n_{\K} \geq 2.
\end{equation}
If $n_{\K}\geq 5$ and $|\Delta_{\K}|$ is sufficiently large, then Louboutin refined \eqref{eq:LouboutinKappa} in \cite{Louboutin15}. We want our results to hold for all number fields, so we favour \eqref{eq:LouboutinKappa} moving forward, although one could use Louboutin's refinement to obtain small improvements to our eventual result when $|\Delta_{\K}|$ is large. Further, if one assumes the Generalised Riemann Hypothesis and $\zeta_{\K}/\zeta$ is entire, then the author and Garcia have established even better explicit bounds for $\kappa_{\K}$ in \cite{GarciaLeeArtinL}. Next, we will prove Lemma \ref{lem:zeta_K_0_absolute_upper_bound_helpful}, which is an observation on the value of $\zeta_{\K}(s)$ at $s=0$. 

\begin{lemma}\label{lem:zeta_K_0_absolute_upper_bound_helpful}
We have $|\zeta_{\K}(0)| \leq A_{\K} \kappa_{\K}$.
\end{lemma}

\begin{proof}
Recall that $\zeta_{\K}(s) = 0$ as long as $r > 0$. If $r=0$, then $\K = \Q$ satisfying $(r_1,r_2) = (1,0)$ \textit{or} $\K$ is an imaginary quadratic field satisfying $(r_1,r_2) = (0,1)$. In the former case, $\zeta_{\K}(0) = -1/2$. In the latter case, \eqref{eqn:functional_equation} and the functional equation for $\Gamma$, $s\Gamma(s) = \Gamma(s+1)$ provide
\begin{align}
    \zeta_{\K}(0) 
    = \lim_{s\to1^{+}} \zeta_{\K}(1-s)
    &= \lim_{s\to1^{+}} \frac{\xi_{\K}(s)}{{A_{\K}}^{1-s} \Gamma\left(1-s\right)^{r_2}}\nonumber\\
    &= \lim_{s\to1^{+}} {A_{\K}}^{2s-1}\frac{\Gamma\left(s\right)}{\Gamma\left(1-s\right)} \zeta_{\K}(s)\nonumber\\
    &= \lim_{s\to1^{+}} {A_{\K}}^{2s-1} (s-1) \frac{\Gamma\left(s-1\right)}{\Gamma\left(1-s\right)} \zeta_{\K}(s)\nonumber\\
    &= - \lim_{s\to1^{+}} {A_{\K}}^{2s-1} (s-1) \zeta_{\K}(s) && \left(\text{since } \lim_{s\to1^{+}} \frac{\Gamma\left(s-1\right)}{\Gamma\left(1-s\right)} = -1\right)\nonumber\\
    &= - A_{\K} \kappa_{\K}. \nonumber 
\end{align}
The result follows naturally.
\end{proof}

Finally, Rademacher used the Phr\'{a}gmen--Lindel\"{o}f theorem in \cite[Thm.~2]{Rademacher} to prove Theorem \ref{thm:RademacherThm4} (below), which re-states \cite[Thm.~4]{Rademacher}.

\begin{theorem}[Rademacher]\label{thm:RademacherThm4}
Suppose that $0 < \eta \leq \tfrac{1}{2}$ and $-\eta \leq \sigma \leq 1 + \eta$, then
\begin{equation*}
    |\zeta_{\K}(s)| \leq 3 \left|\frac{1+s}{1-s}\right| \left(|\Delta_{\K}| \left(\frac{|1+s|}{2\pi}\right)^{n_{\K}}\right)^{\frac{1 + \eta - \sigma}{2}} \zeta(1+\eta)^{n_{\K}},
\end{equation*}
in which $\zeta(s)$ is the Riemann zeta-function.
\end{theorem}

\begin{remark}
In Theorem \ref{thm:RademacherThm4}, the constant $3$ can be replaced by $(1+\eta)/(1-\eta)$. Therefore, there is a small refinement available if one chooses $\eta$ closer to $0$. Alas, we will eventually use Theorem \ref{thm:RademacherThm4} with $\eta = 0.45$, so there is no significant gain to replace the $3$ in this paper. 
\end{remark}

\subsection{The minimum discriminant}

\begin{table}[]
    \centering
    \begin{tabular}{r|lllllll}
        $n_{\K}$ & 2 & 3 & 4 & 5 & 6 & 7 & $\geq 8$\\
        \hline 
        $\Delta_0(n_{\K})$ & 3 & 23 & 117 & 1\,607 & 9\,747 & 184\,607 & $\mathfrak{m}(n_{\K})$ \\
    \end{tabular}
    \caption{Admissible choices for $\Delta_0(n_{\K})$, given a selection of $n_{\K}$, where $\mathfrak{m}(n_{\K})$ is defined in \eqref{eqn:mathfrakmdefined}.}
    \label{tab:SmallestDiscVals}
\end{table}

The invariant $|\Delta_{\K}|$ will play an important role later, because we will need lower bounds for $|\Delta_{\K}|$ to establish several bounds. Let $\Delta_0(n_{\K})$ denote an admissible lower bound for $|\Delta_{\K}|$, in which $\K$ is any number field with degree $n_{\K}$. We restrict our attention to $n_{\K}\geq 2$ without loss of generality, because $\Delta_{\K} = 1$ when $n_{\K} = 1$. For $n_{\K}\geq 2$, Minkowski's well-known bound is
\begin{equation}\label{eqn:mathfrakmdefined}
    |\Delta_{\K}| \geq \left(\frac{\pi}{4}\right)^{n_{\K}} \frac{n_{\K}^{2 n_{\K}}}{(n_{\K}!)^2} := \mathfrak{m}(n_{\K}).
\end{equation}
If $n_{\K} \in [2,7]$, then we know optimum values for these $\Delta_0(n_{\K})$; these values are found by searching the LMDFB database \cite{LMFDB} (which is known to be complete for these degrees). All this information is presented in Table \ref{tab:SmallestDiscVals}.

\begin{remark}
If optimum values for all $\Delta_0(n_{\K})$ such that $n_{\K}\geq 2$ were known, then it is reasonable to think that we should also know $\Delta_0(n_{\K}+1) \geq \Delta_0(n_{\K})$ for each $n_{\K}$. Even though the evidence suggests this is true, we cannot take it for granted. In fact, if it \textit{is} true, then for $n_{\K}\geq 8$ we have $\Delta_0(n_{\K}) = \max\{184\,607, \mathfrak{m}(n_{\K})\}$, so $\Delta_0(n_{\K})$ would improve numerically for $n_{\K}\in \{8,9\}$ and minor refinements in our proof of Theorem \ref{thm:Sunley_improved} would be available.
\end{remark}

\section{Proof of Theorem \ref{thm:Sunley_improved}}\label{sec:MainTheoremProof}

Let $0 < \eta < \tfrac{1}{2}$ be a parameter to be chosen; we ensure $\eta < \tfrac{1}{2}$ so that we can bound $|\zeta_{\K}(-\eta + it)|$ using Theorem \ref{thm:RademacherThm4} later. Recall that $f(s)$ was defined in \eqref{eqn:LandauSatz156} and suppose that
\begin{align*}
    \Omega(x) &= \int_{0}^{x} dy_1 \int_{0}^{y_1} dy_2 \cdots \int_{0}^{y_{n_{\K}-1}} I_{\K}(y_{n_{\K}}) dy_{n_{\K}},\\
    \Upsilon_{\eta}(x) &= \frac{1}{2\pi i} \int_{-\eta - i\infty}^{-\eta + i\infty} \frac{x^{s+n_{\K}}}{s(s+1)\cdots(s+n_{\K})} \zeta_{\K}(s)\,ds,\\
    K_{\eta}(x) &= \frac{1}{2\pi i} \int_{-\eta - i\infty}^{-\eta + i\infty} \frac{x^{s+n_{\K}}}{s(s+1)\cdots(s+n_{\K})} f(s)^{-1}\,ds.
\end{align*}
Write $K(x) = K_{\eta}(x)$ when the choice of $\eta$ does not need to be specified and let $K^{(k)}(w)$ denote the $k$th formal derivative of $K(w)$. As part of his proof of \cite[Satz 207]{LandauEinfuehrung}, Landau shows
\begin{equation*}
    K^{(k)}(w)
    = \frac{1}{2\pi i} \int_{-\tfrac{1}{2} + \tfrac{1}{n_{\K}}\left(k-\tfrac{1}{2}\right) - i\infty}^{-\tfrac{1}{2} + \tfrac{1}{n_{\K}}\left(k-\tfrac{1}{2}\right) + i\infty} \frac{w^{s+n_{\K}-k}}{s \cdots (s+n_{\K}-k)} f(s)^{-1}\,ds .
\end{equation*}
Take $k=0$ and $k=n_{\K}$ respectively to obtain
\begin{align}
    K(w) &= \frac{1}{2\pi i} \int_{-\tfrac{1}{2} - \tfrac{1}{2n_{\K}} - i\infty}^{-\tfrac{1}{2} - \tfrac{1}{2n_{\K}} + i\infty} \frac{w^{s+n_{\K}}}{s \cdots (s+n_{\K})} f(s)^{-1}\,ds,\label{eqn:Kw_defn}\\
    K^{(n_{\K})}(w) &= \frac{1}{2\pi i} \int_{\tfrac{1}{2} - \tfrac{1}{2n_{\K}} - i\infty}^{\tfrac{1}{2} - \tfrac{1}{2n_{\K}} + i\infty} \frac{w^{s}}{s} f(s)^{-1}\,ds.\label{eqn:Knw_defn}
\end{align}

Our approach to prove Theorem \ref{thm:Sunley_improved} follows, and we note that it is the same approach Landau (and Sunley) used to establish \eqref{eqn:Landau} (resp.\ \eqref{eqn:Sunley}) in \cite{LandauEinfuehrung} (resp.\ \cite{SunleyThesis}) with an extra parameter $\eta$. First, one can manipulate the definition of $\Omega(x)$ and move the line of integration to see that
\begin{equation}\label{eqn:UpsilonOmegaRelationship}
    \Omega(x) = \frac{1}{2\pi i} \int_{2 - i\infty}^{2 + i\infty} \frac{x^{s+n_{\K}}}{s(s+1)\cdots(s+n_{\K})} \zeta_{\K}(s)\,ds
    = \Upsilon_{\eta}(x) + \frac{x^{1+n_{\K}}}{(n_{\K}+1)!} \kappa_{\K} + \frac{x^{n_{\K}}}{n_{\K}!} \zeta_{\K}(0),
\end{equation}
Apply \eqref{eqn:DZF_new_def} and $\zeta_{\K}(s) = f(s)^{-1} \zeta_{\K}(1-s)$ in the definition of $\Upsilon_{\eta}(x)$ to obtain
\begin{equation}\label{eqn:observation_Oopsilon}
    \Upsilon_{\eta}(x)
    = \frac{1}{2\pi i} \sum_{m=1}^\infty \frac{\tau_m}{m^{n_{\K}+1}} K_{\eta}(mx).
\end{equation}
Finally, suppose $x > 1$, $z = x^{1 - \tfrac{2}{n_{\K}+1}}$, and
\begin{equation}\label{eqn:ImportantOmegaEstimate2}
    \Xi_z g(w) = \sum_{\ell=0}^{n_{\K}} (-1)^{n_{\K}-\ell} {n_{\K} \choose \ell} g(x + \ell z),
    \quad\text{so}\quad
    z^{n_{\K}} I_{\K}(x) \leq \Xi_z\Omega(x) \leq z^{n_{\K}} I_{\K}(x + n_{\K}z).
\end{equation}
The inequality in \eqref{eqn:ImportantOmegaEstimate2} is the key which unlocks the desired estimate; all we need is an explicit formula for $\Xi_z\Omega(x)$. To this end, we need some explicit ingredients which are given in Section \ref{ssec:ingredients}. Using these ingredients, we obtain an explicit formula for $\Xi_z\Upsilon_{\eta}(x)$ in Section \ref{ssec:xplicitformula}. We complete the proof of Theorem \ref{thm:Sunley_improved} using that explicit formula in Section \ref{ssec:completion}. 

\subsection{Important ingredients}\label{ssec:ingredients}

Recall from \eqref{eq:LouboutinKappa} that
\begin{equation*}
    \kappa_{\K} \leq \alpha(n_{\K})\, (\log{|\Delta_{\K}|})^{n_{\K}-1},
    \qquad\text{where}\qquad
    \alpha(n_{\K}) = \left(\frac{e}{2(n_{\K} - 1)}\right)^{n_{\K} - 1}.
\end{equation*}
We introduce the notation $\alpha(n_{\K})$ for convenience, because it will appear often. The two following theorems are the most important ingredients we require. 

\begin{theorem}\label{thm:Magdeline}
If $x > |\Delta_{\K}|$, then
\begin{equation*}
    I_{\K}(x) \leq Q_1(n_{\K})\, (\log{|\Delta_{\K}|})^{n_{\K}-1} x,
    \quad\text{where}\quad
    Q_1(n_{\K}) 
    = 2.74 \frac{n_{\K} (n_{\K}-1)^{-1}}{(\log{\mathfrak{m}(n_{\K})})^{n_{\K}-1}} n_{\K}!,
\end{equation*}
in which $\mathfrak{m}(n_{\K})$ was defined in \eqref{eqn:mathfrakmdefined}.
\end{theorem}

\begin{theorem}\label{thm:Gertrude}
If $g_{\K}^{\pm} = \frac{5}{8} + \frac{\pi}{2} \pm \frac{1}{n_{\K}} + \frac{3}{8 n_{\K}^2}$ and $w>0$, then 
\begin{align*}
    |K(w)| &\leq Q_2(n_{\K}) |\Delta_{\K}|^{1 + \frac{1}{2\,n_{\K}}} w^{n_{\K}-\frac{1}{2}-\frac{1}{2 n_{\K}}},\\
    |K^{(n_{\K})}(w)| &\leq Q_3(n_{\K}) |\Delta_{\K}|^{\frac{1}{2\,n_{\K}}} w^{\frac{1}{2}-\frac{1}{2 n_{\K}}},
\end{align*}
in which
\begin{align}
    Q_2(n_{\K}) &= 2^{\frac{3}{2}} \pi^{-\frac{n_{\K}}{2}-\frac{3}{2}} e^{4 n_{\K} + \frac{1}{442 n_{\K}} - \frac{1}{2}} \sqrt{g_{\K}^{+}\, n_{\K}} , \nonumber\\ 
    Q_3(n_{\K}) &= 2^{n_{\K} + \frac{3}{2}} \pi^{\frac{n_{\K}}{2}-\frac{3}{2}} e^{n_{\K}\left(1 + \frac{4 n_{\K}}{n_{\K}-1} + \frac{1}{500 n_{\K}^6}\right) - \frac{1}{2}} \sqrt{g_{\K}^{-}\, n_{\K}}. \nonumber 
\end{align}
Note that $g_{\K}^{+} < 2.78955$ and $g_{\K}^{-} < 2.19580$, where the former maximum occurs at $n_{\K}=2$ and the latter maximum is achieved as $n_{\K}\to\infty$.
\end{theorem}

\begin{table}[]
    \centering
    \begin{tabular}{lllllll}
        $n_{\K}$ & $a_7$ & $Q_1(n_{\K})$ & $a_{12}$ & $Q_2(n_{\K})$ & $a_{13}$ & $Q_3(n_{\K})$ \\
        \hline
        2 & $5.1200\cdot 10^{4}$ & $1.2135\cdot 10^{1}$ & $7.8776\cdot 10^{2}$ & $6.9128\cdot 10^{2}$ & $1.6478\cdot 10^{7}$ & $4.8095\cdot 10^{8}$ \\
        3 & $1.2288\cdot 10^{7}$ & $4.7294\cdot 10^{0}$ & $9.1256\cdot 10^{4}$ & $2.5027\cdot 10^{4}$ & $6.5833\cdot 10^{10}$ & $4.3260\cdot 10^{10}$ \\
        4 & $3.9322\cdot 10^{9}$ & $1.6390\cdot 10^{0}$ & $1.2444\cdot 10^{7}$ & $8.7225\cdot 10^{5}$ & $2.9191\cdot 10^{14}$ & $1.3721\cdot 10^{13}$ \\
        5 & $1.5729\cdot 10^{12}$ & $5.1636\cdot 10^{-1}$ & $1.9832\cdot 10^{9}$ & $2.9679\cdot 10^{7}$ & $1.4646\cdot 10^{18}$ & $5.8439\cdot 10^{15}$ \\
        10 & $1.5586\cdot 10^{26}$ & $7.7874\cdot 10^{-4}$ & $1.2770\cdot 10^{21}$ & $1.1366\cdot 10^{15}$ & $2.3614\cdot 10^{37}$ & $1.9530\cdot 10^{29}$ \\
        20 & $1.1220\cdot 10^{57}$ & $4.9142\cdot 10^{-10}$ & $1.9169\cdot 10^{47}$ & $1.2220\cdot 10^{30}$ & $1.5714\cdot 10^{78}$ & $3.5911\cdot 10^{56}$ \\
        30 & $1.3135\cdot 10^{90}$ & $1.7774\cdot 10^{-16}$ & $3.4997\cdot 10^{75}$ & $1.1467\cdot 10^{45}$ & $1.1015\cdot 10^{121}$ & $6.6704\cdot 10^{83}$ \\
        40 & $4.3381\cdot 10^{124}$ & $5.0964\cdot 10^{-23}$ & $1.6013\cdot 10^{105}$ & $1.0165\cdot 10^{60}$ & $1.8243\cdot 10^{165}$ & $1.2102\cdot 10^{111}$ \\
        50 & $1.7363\cdot 10^{160}$ & $1.2854\cdot 10^{-29}$ & $8.3211\cdot 10^{135}$ & $8.7319\cdot 10^{74}$ & $3.3228\cdot 10^{210}$ & $2.1552\cdot 10^{138}$
    \end{tabular}
    \caption{Comparison between values of $a_7$, $a_{12}$, $a_{13}$, $Q_1(n_{\K})$, $Q_2(n_{\K})$, and $Q_3(n_{\K})$ for several choices of $n_{\K}$.}
    \label{tab:comparison_MagdelineGertrude}
\end{table}

Proofs of Theorem \ref{thm:Magdeline} and Theorem \ref{thm:Gertrude} are deferred until Sections \ref{sec:IKxUpperBound} and \ref{sec:KwKnwUpperBounds} respectively, because they are long and technical. In \cite[Thm.~3.2.4]{SunleyThesis}, Sunley proved Theorem \ref{thm:Gertrude} with $a_{12}$ and $a_{13}$ in place of $Q_2(n_{\K})$ and $Q_3(n_{\K})$, where
\begin{align*}
    a_{12} = e^{\frac{7}{2}n_{\K}} \pi^{-\frac{n_{\K}+1}{2}} 2^{\frac{1}{2}-n_{\K}} n_{\K}^{n_{\K}+\frac{3}{2}}
    \quad\text{and}\quad
    a_{13} = e^{5\,n_{\K} + \frac{1}{2}} 2^4 \pi^{\frac{n_{\K}-1}{2}} n_{\K}^{n_{\K}+2}.
\end{align*}
Moreover, she shows in \cite[Thm.~3.1.6]{SunleyThesis} that $a_7 = 2^{4n_{\K}+2} 5^{n_{\K}} {n_{\K}}!$ is admissible in
\begin{equation}\label{eqn:SunleyThm3.1.6}
    I_{\K}(x) \leq a_{\K}(x) (\log{|\Delta_{\K}|})^{n_{\K}-1} x,
    \,\text{ where }\,
    a_{\K}(x) =
    \begin{cases}
    1 & \text{for $0 \leq x < 2$} ,\\[5pt]
    n_{\K} \binom{n_{\K}-1}{\left[ (n_{\K}-1)/2 \right]} & \text{for $2 \leq x \leq |\Delta_{\K}|$} ,\\[5pt]
    a_7 & \text{for $x > |\Delta_{\K}|$} .
    \end{cases}
\end{equation}
\red{Now, let $d_k(m)$ denote the number of
ways of writing $m$ as a product of $k$ integers and note that
\begin{equation}\label{eqn:heswall}
    I_{\K}(x) 
    \leq \#\{ e_1,e_2,\ldots,e_{n_{\K}}\in\Z : e_1e_2\cdots e_{n_{\K}} \leq x \}
    = \sum_{m\leq x} d_{n_{\K}}(m) .
\end{equation}
For any integer $k\geq 1$ and any real $x\geq 1$, Nicolas and Tenenbaum (see \cite[p.~2]{BordellesAO}) also proved that 
\begin{equation*}
    \sum_{m\leq x} d_{k}(m) \leq \frac{x}{(k-1)!} (\log{x} + k-1)^{k-1} ;
\end{equation*}
apply this in \eqref{eqn:heswall} and multiply by one to see that if $2\leq x \leq |\Delta_{\K}|$, then
\begin{equation}\label{eqn:heswall_applied}
    I_{\K}(x)  
    \leq \frac{1}{(n_{\K}-1)!} \left(1 + \frac{n_{\K}-1}{\log{|\Delta_{\K}|}}\right)^{n_{\K}-1} (\log{|\Delta_{\K}|})^{n_{\K}-1} x .
\end{equation}
Therefore, Theorem \ref{thm:Magdeline}, \eqref{eqn:mathfrakmdefined}, and \eqref{eqn:heswall_applied} enable} us to refine \eqref{eqn:SunleyThm3.1.6} into
\begin{equation}\label{eqn:Sunley3.1.6_REFINED}
    I_{\K}(x) \leq b_{\K}(x) (\log{|\Delta_{\K}|})^{n_{\K} - 1}x,
    \quad\text{where}\quad
    b_{\K}(x) =
    \begin{cases}
    1 & \text{for $0 \leq x < 2$} ,\\[5pt]
    \red{\frac{\left(1 + \frac{n_{\K}-1}{\log{\mathfrak{m}(n_{\K})}}\right)^{n_{\K}-1}}{(n_{\K}-1)!}} & \text{for $2 \leq x \leq |\Delta_{\K}|$} ,\\[5pt]
    Q_1(n_{\K}) & \text{for $x > |\Delta_{\K}|$} .
    \end{cases}
\end{equation}
\red{We can also see that $b_{\K}(x) \leq Q_1(n_{\K})$ for all $x\geq 2$.} Now, \eqref{eqn:Sunley3.1.6_REFINED} enables the following explicit version of \cite[Satz 203]{LandauEinfuehrung}.

\begin{lemma}\label{lem:Sunley3.3.1_refined}
If $0 < \theta < 1$ and $x\geq 1$, then
\begin{equation*}
    \sum_{m=1}^{\lfloor x\rfloor} \frac{\tau_m}{m^\theta} < 2\,\red{Q_1(n_{\K})} \left(\frac{1}{1-\theta} + 1\right) (\log{|\Delta_{\K}|})^{n_{\K}-1} x^{1-\theta}.
\end{equation*}
If $\theta > 1$ and $x\geq 1$, then
\begin{equation*}
    \sum_{m=\lfloor x\rfloor+1}^{\infty} \frac{\tau_m}{m^\theta} < \red{Q_1(n_{\K})} \frac{\theta}{\theta - 1} 2^{\theta - 1} (\log{|\Delta_{\K}|})^{n_{\K}-1}  x^{1-\theta}.
\end{equation*}
If $x\geq 2$, then
\begin{equation*}
    \sum_{m=1}^{\lfloor x\rfloor} \frac{\tau_m}{m^\theta} < \red{Q_1(n_{\K})} \left(1 + \frac{2}{\log{2}}\right) (\log{|\Delta_{\K}|})^{n_{\K}-1} \log{x}.
\end{equation*}
\end{lemma}

\begin{proof}
Repeat the proof of \cite[Lem.~3.3.1]{SunleyThesis}, replacing any occurrence of $a_7$ with $Q_1(n_{\K})$.
\end{proof}

Next, an important consequence of Theorem \ref{thm:Gertrude} is the following lemma, which is an explicit version of \cite[Satz 208]{LandauEinfuehrung} that refines \cite[Lem.~3.3.2]{SunleyThesis}.

\begin{lemma}\label{lem:SunleyProp3.3.2_refined}
For $w>1$ and $0 < z < w$, we have
\begin{align*}
    |\Xi_z K(w)| &\leq Q_4(n_{\K})\, |\Delta_{\K}|^{1 + \frac{1}{2\,n_{\K}}} w^{n_{\K}-\frac{1}{2}-\frac{1}{2 n_{\K}}}, \\
    |\Xi_z K(w)| &\leq Q_5(n_{\K})\, z^{n_{\K}} |\Delta_{\K}|^{\frac{1}{2\,n_{\K}}} w^{\frac{1}{2}-\frac{1}{2 n_{\K}}}. 
\end{align*}
in which $Q_4(n_{\K}) = Q_2(n_{\K}) 2^{n_{\K}} (n_{\K}+1)^{n_{\K}-\frac{1}{2}-\frac{1}{2\,n_{\K}}}$ and $Q_5(n_{\K}) = Q_3(n_{\K}) (n_{\K}+1)^{\frac{1}{2} - \frac{1}{2 n_{\K}}}$.
\end{lemma}

\begin{proof}
It follows from Theorem \ref{thm:Gertrude} and $z<w$ that
\begin{align*}
    |\Xi_z K(w)|
    \leq \sum_{\ell=0}^{n_{\K}} {n_{\K} \choose \ell} K(w + \ell z)
    &\leq Q_2(n_{\K}) |\Delta_{\K}|^{1 + \frac{1}{2\,n_{\K}}} \sum_{\ell=0}^{n_{\K}} {n_{\K} \choose \ell} (w + \ell z)^{n_{\K}-\frac{1}{2}-\frac{1}{2\,n_{\K}}}\\
    &< Q_2(n_{\K}) |\Delta_{\K}|^{1 + \frac{1}{2\,n_{\K}}} w^{n_{\K}-\frac{1}{2}-\frac{1}{2\,n_{\K}}} \sum_{\ell=0}^{n_{\K}} {n_{\K} \choose \ell} (\ell + 1)^{n_{\K}-\frac{1}{2}-\frac{1}{2\,n_{\K}}}.
\end{align*}
Now, we observe $(\ell + 1)^{n_{\K}-\frac{1}{2}-\frac{1}{2\,n_{\K}}} \leq (n_{\K} + 1)^{n_{\K}-\frac{1}{2}-\frac{1}{2\,n_{\K}}}$ and import the following standard combinatorial identity:  
\begin{equation}\label{eqn:NIST123}
    \sum_{k=0}^{n_{\K}} \binom{n_{\K}}{k} = 2^{n_{\K}} .
\end{equation}
Apply these to see
\begin{align*}
    |\Xi_z K(w)|
    &\leq Q_2(n_{\K}) (n_{\K}+1)^{n_{\K}-\frac{1}{2}-\frac{1}{2\,n_{\K}}} |\Delta_{\K}|^{1 + \frac{1}{2\,n_{\K}}} w^{n_{\K}-\frac{1}{2}-\frac{1}{2\,n_{\K}}}  \sum_{\ell=0}^{n_{\K}} {n_{\K} \choose \ell}\\
    &= Q_2(n_{\K}) 2^{n_{\K}} (n_{\K}+1)^{n_{\K}-\frac{1}{2}-\frac{1}{2\,n_{\K}}} |\Delta_{\K}|^{1 + \frac{1}{2\,n_{\K}}} w^{n_{\K}-\frac{1}{2}-\frac{1}{2\,n_{\K}}}.
\end{align*}
For the second statement, we use the definition
\begin{equation*}
    \Xi_z K(w) = \int_w^{w+v}dy_1 \cdots \int_{y_{n_{\K}-1}}^{y_{n_{\K}-1}+v} K^{(n_{\K})}(y_{n_{\K}})\,dy_{n_{\K}}.
\end{equation*}
Inserting the upper bound for $K^{(n_{\K})}(w)$ from Theorem \ref{thm:Gertrude} into this definition, we obtain the result by following the steps laid out on \cite[p.~52]{SunleyThesis}.
\end{proof}

Finally, we establish two convenient observations in the following result.

\begin{lemma}\label{lem:MissingDetailsFromLandau}
We have
\begin{equation*}
    \Xi_z(x^{n_{\K}}) = n_{\K}! z^{n_{\K}}
    \qquad\text{and}\qquad
    \Xi_z(x^{n_{\K}+1}) = (n_{\K}+1)! z^{n_{\K}} x + \frac{n_{\K}}{2} (n_{\K}+1)! z^{n_{\K}+1}.
\end{equation*}
\end{lemma}

\begin{proof}
In \cite[Lem.~3.3.3]{SunleyThesis} and \cite[Lem.~3.3.4]{SunleyThesis}, Sunley provides the equality
\begin{equation}\label{eqn:Sunley_lemmas}
    \sum_{\ell=0}^{n_{\K}} (-1)^{n_{\K}-\ell} \binom{n_{\K}}{\ell} \ell^{\alpha} =
    \begin{cases}
        0 & \text{if }0 < \alpha < n_{\K},\\
        n_{\K}! & \text{if }\alpha = n_{\K},\\
        \frac{n_{\K}}{2}(n_{\K}+1)! & \text{if }\alpha = n_{\K}+1,\\
    \end{cases}
\end{equation}
Using \eqref{eqn:Sunley_lemmas} and the binomial theorem, we see that
\begin{align*}
    \Xi_z(x^{n_{\K}})
    = \sum_{\ell=0}^{n_{\K}} (-1)^{n_{\K}-\ell} \binom{n_{\K}}{\ell} (x + \ell z) ^{n_{\K}}
    &= \sum_{\ell=0}^{n_{\K}} (-1)^{n_{\K}-\ell} \binom{n_{\K}}{\ell} \sum_{m=0}^{n_{\K}} \binom{n_{\K}}{m} x^{m} (\ell z) ^{n_{\K}-m}\\
    &= z^{n_{\K}} \sum_{\ell=0}^{n_{\K}} (-1)^{n_{\K}-\ell} \binom{n_{\K}}{\ell} \ell^{n_{\K}}
    = n_{\K}! z^{n_{\K}}.
\end{align*}
Following similar logic, we also see that
\begin{align*}
    \Xi_z(x^{n_{\K}+1})
    &= \sum_{\ell=0}^{n_{\K}} (-1)^{n_{\K}-\ell} \binom{n_{\K}}{\ell} \left(\binom{n_{\K}+1}{n_{\K}} \ell^{n_{\K}} z^{n_{\K}} x + \ell^{n_{\K}+1} z^{n_{\K}+1}\right)\\
    &= (n_{\K}+1)! z^{n_{\K}} x + \frac{n_{\K}}{2} (n_{\K}+1)! z^{n_{\K}+1}. \qedhere
\end{align*}
\end{proof}

\subsection{Explicit formula for $\Xi_z\Omega(x)$}\label{ssec:xplicitformula}

Using the ingredients in the preceding section, we are able to prove the following explicit formula for $\Xi_z\Omega(x)$. Note that the lower bound for the range of $x$ for which our formula holds has been chosen to match a choice in the next section.

\begin{theorem}\label{thm:explicitformulaDeltaOmega}
Suppose that $x > |\Delta_{\K}|^{\frac{1}{2}}$ and $z = |\Delta_{\K}|^{\frac{1}{n_{\K}+1}} x^{\frac{n_{\K}-1}{n_{\K}+1}}$. We have
\begin{equation*}
    \left|\Xi_z\Omega(x) - \kappa_{\K} x z^{n_{\K}}\right| 
    \leq Q_6(x, n_{\K}) (\log{|\Delta_{\K}|})^{n_{\K}-1} z^{n_{\K}+1},
\end{equation*}
such that
\begin{align*}
    &Q_6(x,n_{\K}) \\
    &\quad= \alpha(n_{\K})\left(\frac{n_{\K}}{2} + |\Delta_{\K}|^{\frac{1}{2}\left(\frac{n_{\K}-1}{n_{\K}+1}\right)} x^{-\frac{n_{\K}-1}{n_{\K}+1}} \pi^{-\frac{n_{\K}}{2}}\right) + \frac{Q_1(n_{\K})}{\pi} \left(\frac{6 n_{\K}-2}{n_{\K}-1}\right)\cdot
    \begin{cases}
        Q_5(n_{\K}) &\text{if }n_{\K}\leq 13,\\
        Q_4(n_{\K}) &\text{if }n_{\K}> 13.
    \end{cases}
\end{align*}
\end{theorem}

\begin{proof}
Observe \eqref{eqn:observation_Oopsilon} and split the sum therein at $\lfloor z\rfloor$ to observe
\begin{align*}
    \Xi_{z}\Upsilon_{\eta}(x)
    = \frac{1}{2\pi i} \sum_{m=1}^\infty \frac{\tau_m}{m^{n_{\K}+1}} \Xi_{mz} K_{\eta}(mx)
    = \frac{1}{2\pi i} \left(\sum_{m=1}^{\lfloor z\rfloor} + \sum_{m=\lfloor z\rfloor+1}^\infty\right) \frac{\tau_m}{m^{n_{\K}+1}} \Xi_{mz} K_{\eta}(mx).
\end{align*}
It follows that
\begin{align*}
    |\Xi_{z}\Upsilon_{\eta}(x)| 
    &\leq \frac{1}{2\pi} \left(\sum_{m=1}^{\lfloor z\rfloor} + \sum_{m=\lfloor z\rfloor+1}^\infty\right) \frac{\tau_m}{m^{n_{\K}+1}} |\Xi_{mz} K_{\eta}(mx)|.
\end{align*}
The choice to split the sum at $\lfloor z\rfloor$ was made so that we can implement the convenient observations
\begin{equation*}
    |\Delta_{\K}|^{\frac{1}{2\,n_{\K}}} (xz)^{\frac{1}{2}-\frac{1}{2 n_{\K}}} = z
    \qquad\text{and}\qquad
    |\Delta_{\K}|^{1 + \frac{1}{2 n_{\K}}} x^{n_{\K}} (xz)^{-\frac{1}{2} - \frac{1}{2 n_{\K}}} 
    = z^{n_{\K}+1};
\end{equation*}
these follow using $\tfrac{1}{n_{\K}} - \tfrac{1}{n_{\K}+1} = \tfrac{1}{n_{\K} (n_{\K} + 1)}$. Now, use Lemma \ref{lem:Sunley3.3.1_refined} and Lemma \ref{lem:SunleyProp3.3.2_refined} to see
\begin{align}
    \sum_{m=1}^{\lfloor z\rfloor} \frac{\tau_m}{m^{n_{\K}+1}} |\Xi_{mz} K_{\eta}(mx)|
    &\leq Q_5(n_{\K}) z^{n_{\K}} |\Delta_{\K}|^{\frac{1}{2\,n_{\K}}} x^{\frac{1}{2}-\frac{1}{2 n_{\K}}} \sum_{m=1}^{\lfloor z\rfloor} \frac{\tau_m}{m^{\frac{1}{2}+\frac{1}{2 n_{\K}}}} \nonumber\\
    &\leq Q_1(n_{\K}) Q_5(n_{\K}) (\log{|\Delta_{\K}|})^{n_{\K}-1} z^{n_{\K}} 2\,|\Delta_{\K}|^{\frac{1}{2\,n_{\K}}} \!\!\left(\frac{3 n_{\K} - 1}{n_{\K}-1}\right) (xz)^{\frac{1}{2}-\frac{1}{2 n_{\K}}} \nonumber\\
    &= Q_1(n_{\K}) Q_5(n_{\K}) \left(\frac{6 n_{\K} - 2}{n_{\K}-1}\right) (\log{|\Delta_{\K}|})^{n_{\K}-1} z^{n_{\K}+1}. \label{eqn:Q5UpperBound}
\end{align}
Note that we can apply Lemma \ref{lem:SunleyProp3.3.2_refined}, because $z/x = |\Delta_{\K}|^{\frac{1}{n_{\K}+1}} x^{-\frac{2}{n_{\K}+1}} < 1$, hence $z < x$ for all $x > |\Delta_{\K}|^{1/2}$. Moreover, this is the largest range of $x$ that we could apply Lemma \ref{lem:SunleyProp3.3.2_refined} for under our choice of $z$. Similarly, we have
\begin{align}
    &\sum_{m=\lfloor z\rfloor+1}^\infty \frac{\tau_m}{m^{n_{\K}+1}} |\Xi_{mz} K_{\eta}(mx)| 
    \leq Q_4(n_{\K})|\Delta_{\K}|^{1 + \frac{1}{2 n_{\K}}} x^{n_{\K} -\frac{1}{2} - \frac{1}{2 n_{\K}}} \sum_{m=\lfloor z\rfloor+1}^\infty \frac{\tau_m}{m^{\frac{3}{2} + \frac{1}{2 n_{\K}}}}\nonumber\\
    &\qquad\qquad\qquad\leq Q_1(n_{\K})\,Q_4(n_{\K})\left(\frac{3 n_{\K} + 1}{n_{\K}+1}\right) 2^{\frac{1}{2} + \frac{1}{2 n_{\K}}} |\Delta_{\K}|^{1 + \frac{1}{2 n_{\K}}} (\log{|\Delta_{\K}|})^{n_{\K}-1} x^{n_{\K}} (xz)^{-\frac{1}{2} - \frac{1}{2 n_{\K}}}\nonumber\\
    &\qquad\qquad\qquad= 2^{\frac{1}{2} + \frac{1}{2 n_{\K}}} \left(\frac{3 n_{\K} + 1}{n_{\K}+1}\right) Q_1(n_{\K})\,Q_4(n_{\K}) (\log{|\Delta_{\K}|})^{n_{\K}-1} z^{n_{\K}+1}. \label{eqn:Q4UpperBound}
\end{align}
It follows from \eqref{eqn:Q5UpperBound} and \eqref{eqn:Q4UpperBound} that
\begin{align}
    |\Xi_{z}\Upsilon_{\eta}(x)|
    &\leq \frac{1}{\pi} \left(\frac{6 n_{\K}-2}{n_{\K}-1}\right) Q_1(n_{\K}) \max\left\{Q_4(n_{\K}), Q_5(n_{\K})\right\} (\log{|\Delta_{\K}|})^{n_{\K}-1} z^{n_{\K}+1} \nonumber\\
    &= \frac{Q_1(n_{\K})}{\pi} \left(\frac{6 n_{\K}-2}{n_{\K}-1}\right) (\log{|\Delta_{\K}|})^{n_{\K}-1} z^{n_{\K}+1} \cdot \begin{cases}
         Q_5(n_{\K}) &\text{if }n_{\K} \leq  13,\\
         Q_4(n_{\K}) &\text{if }n_{\K} >  13.
    \end{cases}\label{eqn:ENDLICH}
\end{align}
Finally, use \eqref{eqn:UpsilonOmegaRelationship}, \eqref{eqn:ENDLICH}, and Lemma \ref{lem:MissingDetailsFromLandau} to obtain
\begin{align*}
    \Xi_z\Omega(x)
    &= \Xi_z\Upsilon_{\eta}(x) + \frac{\Xi_z(x^{1+n_{\K}})}{(n_{\K}+1)!} \kappa_{\K} + \frac{\Xi_z(x^{n_{\K}})}{n_{\K}!} \zeta_{\K}(0)\\
    &= \kappa_{\K} x z^{n_{\K}} + \frac{n_{\K}}{2} \kappa_{\K} z^{n_{\K}+1} + \zeta_{\K}(0) z^{n_{\K}} \\
    &\qquad\qquad+ 
    \begin{cases}
        O^*\!\!\left( \frac{1}{\pi} \left(\frac{6 n_{\K}-2}{n_{\K}-1}\right)  Q_1(n_{\K}) Q_5(n_{\K}) (\log{|\Delta_{\K}|})^{n_{\K}-1} z^{n_{\K}+1} \right) &\text{if }n_{\K}\leq 13,\\
        O^*\!\!\left( \frac{1}{\pi} \left(\frac{6 n_{\K}-2}{n_{\K}-1}\right)  Q_1(n_{\K}) Q_4(n_{\K}) (\log{|\Delta_{\K}|})^{n_{\K}-1} z^{n_{\K}+1} \right)&\text{if }n_{\K}> 13.
    \end{cases}
\end{align*}
It follows from \eqref{eq:LouboutinKappa} and Lemma \ref{lem:zeta_K_0_absolute_upper_bound_helpful} that
\begin{align*}
    &|\Xi_z\Omega(x) - \kappa_{\K} x z^{n_{\K}}|\\
    &\leq \frac{n_{\K}}{2} \kappa_{\K} z^{n_{\K}+1} + |\zeta_{\K}(0)| z^{n_{\K}} + 
    \frac{Q_1(n_{\K})}{\pi} \left(\frac{6 n_{\K}-2}{n_{\K}-1}\right)  (\log{|\Delta_{\K}|})^{n_{\K}-1} z^{n_{\K}+1}\cdot
    \begin{cases}
        Q_5(n_{\K}) &\text{if }n_{\K}\leq 13\\
        Q_4(n_{\K}) &\text{if }n_{\K}> 13
    \end{cases}\\
    &\leq Q_6(x, n_{\K}) (\log{|\Delta_{\K}|})^{n_{\K}-1} z^{n_{\K}+1}. \qedhere
\end{align*}
\end{proof}

\subsection{Completing the proof of Theorem \ref{thm:Sunley_improved}}\label{ssec:completion}

We establish the result for $x \leq (n_{\K}+1) |\Delta_{\K}|^{\frac{1}{2}}$ and $x > (n_{\K}+1) |\Delta_{\K}|^{\frac{1}{2}}$ separately. Note that we split this range here, to ensure that the lower range yields a bound of the desired order.

\medskip
\noindent\textsc{Case I:} Suppose that $0 < x \leq (n_{\K}+1) |\Delta_{\K}|^{\tfrac{1}{2}}$.
Use \eqref{eq:LouboutinKappa} and Theorem \ref{thm:Magdeline} to see
\begin{align*}
    |I_{\K}(x) - \kappa_{\K} x|
    &\leq |I_{\K}(x)| + |\kappa_{\K} x|\\
    &\leq \left(Q_1(n_{\K}) +\alpha(n_{\K})\right) (\log{|\Delta_{\K}|})^{n_{\K}-1} x\\
    &= \left(Q_1(n_{\K}) + \alpha(n_{\K})\right) (\log{|\Delta_{\K}|})^{n_{\K}-1} x^{\frac{n_{\K}-1}{n_{\K}+1}} x^{\frac{2}{n_{\K}+1}}\\
    &\leq \left(Q_1(n_{\K}) + \alpha(n_{\K})\right) \left(n_{\K}+1\right)^{\frac{2}{n_{\K}+1}} |\Delta_{\K}|^{\frac{1}{n_{\K}+1}} (\log{|\Delta_{\K}|})^{n_{\K}-1}  x^{\frac{n_{\K}-1}{n_{\K}+1}}.
\end{align*}
Now, the term $\left(Q_1(n_{\K}) + \alpha(n_{\K})\right) \left(n_{\K}+1\right)^{\frac{2}{n_{\K}+1}}$ decreases for $n_{\K} \geq 2$ bound. Moreover, it is clearly majorised by $\Lambda_{\K}(n_{\K})$, because $\Lambda_{\K}(n_{\K})$ increases for $n_{\K}\geq 2$ and its value at $n_{\K}=2$ is less than $\Lambda_{\K}(2)$.

\medskip
\noindent\textsc{Case II:} Suppose that $x > |\Delta_{\K}|^{\frac{1}{2}}$, and $z = |\Delta_{\K}|^{\frac{1}{n_{\K}+1}} x^{\frac{n_{\K}-1}{n_{\K}+1}}$. Use the inequality in \eqref{eqn:ImportantOmegaEstimate2} and Theorem \ref{thm:explicitformulaDeltaOmega} to see
\begin{align*}
    I_{\K}(x)
    \leq z^{-n_{\K}} \Xi_z\Omega(x)
    &\leq \kappa_{\K} x + Q_6(x, n_{\K})\, (\log{|\Delta_{\K}|})^{n_{\K}-1} z\\
    &= \kappa_{\K} x + Q_6(x, n_{\K})\, |\Delta_{\K}|^{\frac{1}{n_{\K}+1}} (\log{|\Delta_{\K}|})^{n_{\K}-1} x^{1 - \frac{2}{n_{\K}+1}}.
\end{align*}
Similarly, 
\begin{align*}
    I_{\K}(x + n_{\K}z)
    \geq z^{-n_{\K}} \Xi_z\Omega(x)
    &\geq \kappa_{\K} x - Q_6(x, n_{\K})\, |\Delta_{\K}|^{\frac{1}{n_{\K}+1}} (\log{|\Delta_{\K}|})^{n_{\K}-1} x^{1 - \frac{2}{n_{\K}+1}}.
\end{align*}
Therefore, write $y = x + n_{\K}z$, then $y > |\Delta_{\K}|^{\frac{1}{2}} + n_{\K} |\Delta_{\K}|^{\frac{1}{n_{\K}+1}} |\Delta_{\K}|^{\frac{1}{2} - \frac{1}{n_{\K}+1}} = (n_{\K} + 1) |\Delta_{\K}|^{\tfrac{1}{2}}$ and
\begin{align*}
    I_{\K}(y)
    &> \kappa_{\K} \left(y - n_{\K} z\right) - Q_6(x, n_{\K})\, |\Delta_{\K}|^{\frac{1}{n_{\K}+1}} (\log{|\Delta_{\K}|})^{n_{\K}-1} y^{1 - \frac{2}{n_{\K}+1}} \\
    &> \kappa_{\K} y - \left(Q_6(x, n_{\K}) + \frac{n_{\K} \kappa_{\K}}{(\log{|\Delta_{\K}|})^{n_{\K}-1}}\right) |\Delta_{\K}|^{\frac{1}{n_{\K}+1}} (\log{|\Delta_{\K}|})^{n_{\K}-1} y^{1 - \frac{2}{n_{\K}+1}}.
\end{align*}
Combine these observations with \eqref{eq:LouboutinKappa} to see that for $x > (n_{\K}+1) |\Delta_{\K}|^{\frac{1}{2}}$,
\begin{align*}
    |I_{\K}(x) - \kappa_{\K} x | 
    &< \left(Q_6(x, n_{\K}) + n_{\K} \alpha(n_{\K})\right) |\Delta_{\K}|^{\frac{1}{n_{\K}+1}} (\log{|\Delta_{\K}|})^{n_{\K}-1} x^{1 - \frac{2}{n_{\K}+1}}.
\end{align*}
All that remains to complete the proof of Theorem \ref{thm:Sunley_improved} is to show that $Q_6(x, n_{\K}) + n_{\K} \alpha(n_{\K})$ is bounded from above by $\Lambda_{\K}(n_{\K})$. That is, for $x > (n_{\K}+1) |\Delta_{\K}|^{\frac{1}{2}}$ we have
\begin{align*}
    &Q_6(x, n_{\K}) + n_{\K} \alpha(n_{\K})\\
    &\quad\leq \alpha(n_{\K}) \left(\frac{3 n_{\K}}{2} + |\Delta_{\K}|^{\frac{1}{2} - \frac{1}{n_{\K}+1}} x^{-\frac{n_{\K}-1}{n_{\K}+1}} \pi^{-\frac{n_{\K}}{2}}\right) + \frac{Q_1(n_{\K})}{\pi} \left(\frac{6 n_{\K}-2}{n_{\K}-1}\right) \cdot
    \begin{cases}
        Q_5(n_{\K}) &\text{if }n_{\K}\leq 13\\
        Q_4(n_{\K}) &\text{if }n_{\K}> 13
    \end{cases} \\
    &\quad< \alpha(n_{\K}) \left(\frac{3 n_{\K}}{2} + (n_{\K}+1)^{-\frac{n_{\K}-1}{n_{\K}+1}} \pi^{-\frac{n_{\K}}{2}}\right) + \frac{Q_1(n_{\K})}{\pi} \left(\frac{6 n_{\K}-2}{n_{\K}-1}\right) \cdot
    \begin{cases}
        Q_5(n_{\K}) &\text{if }n_{\K}\leq 13\\
        Q_4(n_{\K}) &\text{if }n_{\K}> 13
    \end{cases} \\
    &\quad\leq \frac{e^{0.01 n_{\K}^{-1}} Q_1(n_{\K})}{\pi} \left(\frac{6 n_{\K}-2}{n_{\K}-1}\right) \cdot
    \begin{cases}
        Q_5(n_{\K}) &\text{if }n_{\K}\leq 13,\\
        Q_4(n_{\K}) &\text{if }n_{\K}> 13.
    \end{cases}.
\end{align*}
To justify the last inequality for $n_{\K} \leq 13$, we confirmed it using computations. To justify it for $n_{\K} \geq 14$, note that the following ratio is decreasing for $n_{\K} \geq 14$:
\begin{align*}
    \frac{\alpha(n_{\K}) \left(\frac{3 n_{\K}}{2} + (n_{\K}+1)^{-\frac{n_{\K}-1}{n_{\K}+1}} \pi^{-\frac{n_{\K}}{2}}\right)}{\frac{Q_1(n_{\K})}{\pi} \left(\frac{6 n_{\K}-2}{n_{\K}-1}\right) Q_4(n_{\K}) \left(e^{\frac{0.01}{n_{\K}}}-1\right)} .
\end{align*}
It follows that this ratio is majorised by its value at $n_{\K} = 14$ for all $n_{\K} \geq 14$, which is subsequently majorised by one, so 
\begin{align*}
    \alpha(n_{\K}) \left(\frac{3 n_{\K}}{2} + (n_{\K}+1)^{-\frac{n_{\K}-1}{n_{\K}+1}} \pi^{-\frac{n_{\K}}{2}}\right) < \left(e^{\frac{0.01}{n_{\K}}}-1\right) \frac{Q_1(n_{\K})}{\pi} \left(\frac{6 n_{\K}-2}{n_{\K}-1}\right) Q_4(n_{\K}).
\end{align*}
Now, if $n_{\K} \leq 13$, then
\begin{align*}
    &Q_6(x, n_{\K}) + n_{\K} \alpha(n_{\K})\\
    &\leq 2.74 \frac{2^{\frac{5}{2}}\sqrt{g_{\K}^{-}}}{\pi^{\frac{5}{2}} \sqrt{e}} \frac{(3 n_{\K}-1) (n_{\K}+1)^{\frac{1}{2} - \frac{1}{2 n_{\K}}} }{(n_{\K}-1)^2(\log{\mathfrak{m}(n_{\K})})^{n_{\K}-1}} \left(2\sqrt{\pi}\right)^{n_{\K}} e^{n_{\K}\left(1 + \frac{4 n_{\K}}{n_{\K}-1} + \frac{0.01}{n_{\K}^2} + \frac{1}{500 n_{\K}^6}\right)} n_{\K}^{\frac{3}{2}}n_{\K}! \\
    &\leq 0.53741 \frac{(3 n_{\K}-1) (n_{\K}+1)^{\frac{1}{2} - \frac{1}{2 n_{\K}}} }{(n_{\K}-1)^2(\log{\mathfrak{m}(n_{\K})})^{n_{\K}-1}} \left(\frac{5}{8} + \frac{\pi}{2} - \frac{1}{n_{\K}} + \frac{3}{8 n_{\K}^2}\right)^{\frac{1}{2}} e^{n_{\K}\left(2.26552 + \frac{4 n_{\K}}{n_{\K}-1} + \frac{0.01}{n_{\K}^2} + \frac{1}{500 n_{\K}^6}\right)} n_{\K}^{\frac{3}{2}} n_{\K}! ,
\end{align*}
because $2\sqrt{\pi} \leq e^{1.26552}$. Similarly, if $n_{\K} > 13$, then
\begin{align*}
    &Q_6(x, n_{\K}) + n_{\K} \alpha(n_{\K})\\
    &\leq 2.74 \frac{2^{\frac{5}{2}} \sqrt{g_{\K}^{+}}}{\pi^{\frac{5}{2}}\sqrt{e}} \frac{(3 n_{\K}-1) (n_{\K}+1)^{n_{\K}-\frac{1}{2}-\frac{1}{2\,n_{\K}}}}{(n_{\K}-1)^2 (\log{\mathfrak{m}(n_{\K})})^{n_{\K}-1}} \left(\frac{2}{\sqrt{\pi}}\right)^{n_{\K}} e^{4 n_{\K} + \frac{0.01227}{n_{\K}}} n_{\K}^{\frac{3}{2}} n_{\K}! \\
    &\leq 0.53741 \frac{(3 n_{\K}-1) (n_{\K}+1)^{n_{\K}-\frac{1}{2}-\frac{1}{2\,n_{\K}}}}{(n_{\K}-1)^2 (\log{\mathfrak{m}(n_{\K})})^{n_{\K}-1}} \left(\frac{5}{8} + \frac{\pi}{2} + \frac{1}{n_{\K}} + \frac{3}{8 n_{\K}^2}\right)^{\frac{1}{2}} e^{4.12079 n_{\K} + \frac{0.01227}{n_{\K}}} n_{\K}^{3/2} n_{\K}!, 
\end{align*}
because $2/\sqrt{\pi} \leq e^{0.12079}$. Combining these observations establishes Theorem \ref{thm:Sunley_improved}.

\section{Proof of Theorem \ref{thm:Magdeline}}\label{sec:IKxUpperBound}

Suppose $n_{\K} \geq 2$ and $x > |\Delta_{\K}| \geq \Delta_0(n_{\K})$, where $\Delta_0(n_{\K})$ were defined in Table \ref{tab:SmallestDiscVals}. The inequality $|\Delta_{\K}| \geq \Delta_0(n_{\K})$ will be important, because it enables refinements to our computations. Three preliminary estimates are required to prove Theorem \ref{thm:Magdeline}. The first is $\Upsilon_{\eta}(2x) \ll x^{n_{\K}+1}$ with explicit constants.

\begin{lemma}\label{lem:Madge1}
For $0.01\leq \eta\leq 0.45$, $x>|\Delta_{\K}|^{\tfrac{1}{2}}$, and $$v_1(\eta, n_{\K}) = \frac{9.98134\, \zeta(1+\eta)^{n_{\K}} 2^{n_{\K}-\eta} n_{\K}!}{\Delta_0(n_{\K})^{\tfrac{1}{2}} (\log{\Delta_0(n_{\K})})^{n_{\K}-1} (2\pi)^{n_{\K}\left(\tfrac{1}{2} + \eta\right)}},$$ we have
\begin{equation*}
    \frac{n_{\K}!}{x^{n_{\K}}} |\Upsilon_{\eta}(2x)|
    \leq v_1(\eta, n_{\K}) (\log{|\Delta_{\K}|})^{n_{\K}-1} x. 
\end{equation*}
\end{lemma}

\begin{proof}
Theorem \ref{thm:RademacherThm4} implies
\begin{align*}
    |\Upsilon_{\eta}(2x)|
    &= \left|\frac{1}{2\pi i} \int_{-\eta - i\infty}^{-\eta + i\infty} \frac{(2x)^{s+n_{\K}}}{s(s+1)\cdots(s+n_{\K})} \zeta_{\K}(s)\,ds\right|\\
    &\leq \frac{3}{2\pi} (2x)^{n_{\K} - \eta} \left(\frac{|\Delta_{\K}|}{(2\pi)^{n_{\K}}}\right)^{\frac{1}{2} + \eta} \zeta(1+\eta)^{n_{\K}} \int_{-\eta - i\infty}^{-\eta + i\infty} \frac{\left|\frac{1+s}{1-s}\right| |1+s|^{n_{\K}\left(\frac{1}{2} + \eta\right)}}{|s| |s+1| \cdots |s+n_{\K}|}\,ds\\
    &< \frac{3}{2\pi} (2x)^{n_{\K} - \eta} \left(\frac{|\Delta_{\K}|}{(2\pi)^{n_{\K}}}\right)^{\frac{1}{2} + \eta} \zeta(1+\eta)^{n_{\K}} \int_{-\eta - i\infty}^{-\eta + i\infty}
    \frac{|1+s| |1+s|^{n_{\K}\left(\frac{1}{2} + \eta\right)}}{|s| |1-s| |1 + s|^{n_{\K}}}\,ds.
\end{align*}
Herein, we can estimate the integral as follows:
\begin{align*}
    \int_{-\eta - i\infty}^{-\eta + i\infty}
    \frac{|1+s| |1+s|^{n_{\K}\left(\frac{1}{2} + \eta\right)}}{|s| |1-s| |1 + s|^{n_{\K}}}\,ds
    &\leq \int_{-\eta - i\infty}^{-\eta + i\infty}
    \frac{|1+s|}{|s| |1-s|} |1+s|^{n_{\K}\left(\eta -\frac{1}{2}\right)}\,ds\\
    &\leq \int_{-\eta - i\infty}^{-\eta + i\infty}
    \frac{|1+s|}{|s| |1-s|} |1+s|^{-1 + 2\eta}\,ds\\
    &= \int_{-\eta - i\infty}^{-\eta + i\infty}
    \frac{|1+s|^{2\eta}}{|s| |1-s|}\,ds\\
    &= 2 \int_{0}^{\infty}
    \frac{|1-\eta + it|^{2\eta}}{|-\eta + it| |1+\eta - it|}\,dt .
\end{align*}
This integral converges for $0 < \eta < \tfrac{1}{2}$, and computations show that for $0.01 \leq \eta \leq 0.45$,
\begin{equation*}
    \int_{0}^{\infty}
    \frac{|1-\eta + it|^{2\eta}}{|-\eta + it| |1+\eta - it|}\,dt < 10.45243,
\end{equation*}
where the maximum value is obtained at $\eta = 0.45$. Now, use $x > |\Delta_{\K}|$ to see
\begin{align*}
    |\Upsilon_{\eta}(2x)|
    &< 9.98134 \left(\frac{|\Delta_{\K}|}{(2\pi)^{n_{\K}}}\right)^{\frac{1}{2} + \eta} \zeta(1+\eta)^{n_{\K}} (2x)^{n_{\K}-\eta} \nonumber\\ 
    &< 9.98134 \frac{\zeta(1+\eta)^{n_{\K}} 2^{n_{\K}-\eta}}{x^{\tfrac{1}{2}} (\log{|\Delta_{\K}|})^{n_{\K}-1} (2\pi)^{n_{\K}\left(\tfrac{1}{2} + \eta\right)}} \,(\log{|\Delta_{\K}|})^{n_{\K}-1} x^{n_{\K} + 1}\\
    &< 9.98134 \frac{\zeta(1+\eta)^{n_{\K}} 2^{n_{\K}-\eta}}{|\Delta_{\K}|^{\tfrac{1}{2}} (\log{|\Delta_{\K}|})^{n_{\K}-1} (2\pi)^{n_{\K}\left(\tfrac{1}{2} + \eta\right)}} \,(\log{|\Delta_{\K}|})^{n_{\K}-1} x^{n_{\K} + 1}. 
\end{align*}
Observe $|\Delta_{\K}|\geq \Delta_0(n_{\K})$ to reveal the result.
\end{proof}

We also obtain Lemmas \ref{lem:Madge2} and \ref{lem:Madge3} using \eqref{eq:LouboutinKappa}.

\begin{lemma}\label{lem:Madge2}
For $v_2(n_{\K}) = \frac{4}{n_{\K}+1} \left(\frac{e}{n_{\K}-1}\right)^{n_{\K}-1}$, we have
\begin{equation*}
    \frac{2^{n_{\K}+1}}{n_{\K}+1} \kappa_{\K} 
    \leq v_2(n_{\K}) (\log{|\Delta_{\K}|})^{n_{\K}-1}.
\end{equation*}
\end{lemma}

\begin{lemma}\label{lem:Madge3}
For $v_3(n_{\K}) = \frac{2\, \pi^{-\frac{n_{\K}}{2}}}{\sqrt{\Delta_0(n_{\K})}}  \left(\frac{e}{n_{\K}-1}\right)^{n_{\K}-1}$, we have
\begin{equation*}
    2^{n_{\K}} |\zeta_{\K}(0)| \leq v_3(n_{\K}) (\log{|\Delta_{\K}|})^{n_{\K}-1} x.
\end{equation*}
\end{lemma}

\begin{proof}
Use \eqref{eq:LouboutinKappa} in Lemma \ref{lem:zeta_K_0_absolute_upper_bound_helpful} and $x > |\Delta_{\K}| \geq \Delta_0(n_{\K})$.
\end{proof}

Finally, we need the key observation
\begin{align*}
    |\Omega(2x)|
    = \left|\int_{0}^{2x} dy_1 \int_{0}^{y_1} dy_2 \cdots \int_{0}^{y_{n_{\K}-1}} I_{\K}(y_{n_{\K}}) dy_{n_{\K}}\right|
    &\geq \left|\int_{x}^{2x} dy_1 \int_{x}^{y_1} dy_2 \cdots \int_{x}^{y_{n_{\K}-1}} I_{\K}(y_{n_{\K}}) dy_{n_{\K}}\right|\\
    &\geq \frac{x^{n_{\K}}}{n_{\K}!} I_{\K}(x).
\end{align*}
Use this and \eqref{eqn:UpsilonOmegaRelationship} to see
\begin{align}
    I_{\K}(x) 
    \leq \frac{n_{\K}!}{x^{n_{\K}}} |\Omega(2x)|
    &= \frac{n_{\K}!}{x^{n_{\K}}} |\Upsilon_{\eta}(2x)| +  \frac{2^{n_{\K}+1}}{n_{\K}+1} \kappa_{\K} x  + 2^{n_{\K}} |\zeta_{\K}(0)|,\label{eqn:IK_Omega_relationship}
\end{align}
Choose $\eta = 0.45$ and estimate each term in \eqref{eqn:IK_Omega_relationship} using Lemmas \ref{lem:Madge1}, \ref{lem:Madge2}, and \ref{lem:Madge3} to see
\begin{align*}
    I_{\K}(x)
    &\leq \left(v_1(0.45,n_{\K}) + v_2(n_{\K}) + v_3(n_{\K})\right) (\log{|\Delta_{\K}|})^{n_{\K}-1} x \\ 
    &\leq Q_1(n_{\K}) (\log{|\Delta_{\K}|})^{n_{\K}-1} x. 
\end{align*}
This final upper bound has been confirmed by computations for $n_{\K} \leq 7$. Further, to justify the bound for $n_{\K} \geq 8$, note that $\wp(n_{\K})$ decreases in $n_{\K}$ and hence $\wp(n_{\K}) \leq \wp(8) < 2.74$, where
\begin{align*}
    \wp(n_{\K}) = \frac{\left(v_1(0.45,n_{\K}) + v_2(n_{\K}) + v_3(n_{\K})\right) (n_{\K}-1) (\log{\mathfrak{m}(n_{\K})})^{n_{\K}-1}}{n_{\K}! n_{\K}}. 
\end{align*}
A selection of computations is presented in Table \ref{tab:justificationQ1}. \qed

\begin{table}
    \centering
    \begin{tabular}{l|lllll}
        $n_{\K}$ & $v_1(0.45,n_{\K})$ & $v_2(n_{\K})$ & $v_3(n_{\K})$ & $v_1 + v_2 + v_3$ & $Q_1(n_{\K})$\\
        \hline
        2 & $7.49564\cdot 10^{0}$ & $3.62438\cdot 10^{0}$ & $9.99112\cdot 10^{-1}$ & $1.21191\cdot 10^{1}$ & $1.21351\cdot 10^{1}$ \\
        3 & $8.96578\cdot 10^{-1}$ & $1.84726\cdot 10^{0}$ & $1.38347\cdot 10^{-1}$ & $2.88219\cdot 10^{0}$ & $4.72940\cdot 10^{0}$ \\
        4 & $1.43003\cdot 10^{-1}$ & $5.95127\cdot 10^{-1}$ & $1.39366\cdot 10^{-2}$ & $7.52066\cdot 10^{-1}$ & $1.63897\cdot 10^{0}$ \\
        5 & $6.92239\cdot 10^{-3}$ & $1.42183\cdot 10^{-1}$ & $6.07876\cdot 10^{-4}$ & $1.49713\cdot 10^{-1}$ & $5.16356\cdot 10^{-1}$ \\
        6 & $7.57999\cdot 10^{-4}$ & $2.71384\cdot 10^{-2}$ & $3.10290\cdot 10^{-5}$ & $2.79274\cdot 10^{-2}$ & $1.52108\cdot 10^{-1}$ \\
        7 & $2.47647\cdot 10^{-5}$ & $4.32344\cdot 10^{-3}$ & $7.32387\cdot 10^{-7}$ & $4.34894\cdot 10^{-3}$ & $4.26741\cdot 10^{-2}$ \\
        8 & $1.54330\cdot 10^{-4}$ & $5.91824\cdot 10^{-4}$ & $1.72682\cdot 10^{-7}$ & $7.46327\cdot 10^{-4}$ & $1.15397\cdot 10^{-2}$ \\
        9 & $1.78992\cdot 10^{-5}$ & $7.10716\cdot 10^{-5}$ & $5.71699\cdot 10^{-9}$ & $8.89765\cdot 10^{-5}$ & $3.03231\cdot 10^{-3}$ \\
        10 & $2.01009\cdot 10^{-6}$ & $7.60563\cdot 10^{-6}$ & $1.65982\cdot 10^{-10}$ & $9.61588\cdot 10^{-6}$ & $7.78740\cdot 10^{-4}$ \\
        25 & $2.07974\cdot 10^{-21}$ & $3.05551\cdot 10^{-24}$ & $8.67314\cdot 10^{-38}$ & $2.08280\cdot 10^{-21}$ & $3.07619\cdot 10^{-13}$ \\
        50 & $2.63033\cdot 10^{-47}$ & $2.26627\cdot 10^{-63}$ & $3.09337\cdot 10^{-92}$ & $2.63033\cdot 10^{-47}$ & $1.28536\cdot 10^{-29}$ \\
        75 & $1.60305\cdot 10^{-73}$ & $3.43472\cdot 10^{-108}$ & $1.48491\cdot 10^{-152}$ & $1.60305\cdot 10^{-73}$ & $3.02671\cdot 10^{-46}$
    \end{tabular}
    \caption{A comparison of $v_1(0.45, n_{\K}) + v_2(n_{\K}) + v_3(n_{\K})$ against $Q_1(n_{\K})$ for several choices of $n_{\K}$.}
    \label{tab:justificationQ1}
\end{table}

\section{Proof of Theorem \ref{thm:Gertrude}}\label{sec:KwKnwUpperBounds}

To prove Theorem \ref{thm:Gertrude}, we follow similar techniques to those in the proof of \cite[Thm.~3.2.4]{SunleyThesis}, with some modifications. To this end, let 
\begin{equation*}
    w>0,\quad 
    \mu_{\K} = -\tfrac{1}{2} - \tfrac{1}{2n_{\K}},\quad 
    z_1 = \mu_{\K} + it, \quad
    \nu_{\K} = \tfrac{1}{2} - \tfrac{1}{2n_{\K}}, \quad
    \text{and}\quad
    z_2 = \nu_{\K} + it.
\end{equation*}
In what follows, we collect important preliminary results together in Section \ref{sec:POGertrude}, then use those observations to deduce Theorem \ref{thm:Gertrude} in Section \ref{sec:ProofGertrude}.

\subsection{Preliminary observations}\label{sec:POGertrude}

The first step in our approach requires \eqref{eqn:Kw_defn} and \eqref{eqn:Knw_defn}, so we will need to estimate
\begin{align*}
    f(s)^{-1} &= \left(\frac{(2\pi)^s}{2}\right)^{n_{\K}} |\Delta_{\K}|^{\frac{1}{2}-s}\left(\cos{\frac{\pi s}{2}}\right)^{-r_1 - r_2} \left(\sin{\frac{\pi s}{2}}\right)^{-r_2} \Gamma(s)^{-n_{\K}},
\end{align*}
such that $\Re{s}\in\left\{\mu_{\K}, \nu_{\K}\right\}$. To this end, make the definitions $s=\sigma + it$ with $t\geq 0$,
\begin{align*}
    z_s = e^{-\pi t + \pi \sigma i},\quad
    \cos\frac{\pi s}{2} 
    = \frac{e^{\frac{\pi}{2}\left(t-\sigma i\right)}}{2} \left(1 + z_s\right),
    \quad\text{and}\quad
    \sin\frac{\pi s}{2} 
    = \frac{e^{\frac{\pi}{2}\left(t-\sigma i\right)} i}{2} \left(1 - z_s\right).
\end{align*}
Apply these definitions and $n_{\K} = r_1 + 2\,r_2$ to obtain
\begin{align}
    \left(\cos\frac{\pi s}{2}\right)^{-r_1-r_2} \left(\sin\frac{\pi s}{2}\right)^{-r_2}
    &= \frac{2^{n_{\K}}}{i^{r_2}}\, e^{- \frac{n_{\K} \pi}{2} \left(t-\sigma i\right)} \left(\frac{1}{1 + z_s}\right)^{r_1+r_2} \left(\frac{1}{1 - z_s}\right)^{r_2} \nonumber\\
    &= \frac{2^{n_{\K}}}{i^{r_2}}\, e^{- \frac{n_{\K} \pi}{2} \left(t-\sigma i\right)} \left(1 + \alpha_{\K}(t)\right), \label{eqn:defsapplied}  
\end{align}
where $\alpha_{\K}(t) = \left(\frac{1}{1 + z_s}\right)^{r_1+r_2} \left(\frac{1}{1 - z_s}\right)^{r_2} - 1$. Next, equation (3) on \cite[p.~36]{SunleyThesis} tells us
\begin{equation}\label{eqn:GammaInsert}
    \Gamma(s)^{-n_{\K}} = t^{n_{\K} \left(\frac{1}{2} - \sigma\right)} \exp\left(\frac{\pi n_{\K}}{2}t + \frac{\pi n_{\K}}{2} \left(\frac{1}{2} - \sigma\right)i - n_{\K} c - (\log{t} - 1) n_{\K} ti + \beta(\sigma,t)\right),
\end{equation}
such that $|c| \leq \tfrac{1}{2}$ and
\begin{equation}\label{eqn:beta_defn}
    \beta(\sigma,t) = n_{\K} \left(\frac{1}{2} - \sigma\right) \int_0^{\sigma} \frac{du}{u+ti} + n_{\K} \int_0^{\sigma} \frac{u\,du}{u+ti} - n_{\K} \int_0^{\infty} \frac{u - \lfloor u\rfloor - 1/2}{u+s}\,du .
\end{equation}
Note that \eqref{eqn:GammaInsert} is Stirling's formula with an explicit error term. Finally, we will need explicit bounds for
\begin{align*}
    I &= \int_{0}^{\infty} \frac{t^{n_{\K}\left(\frac{1}{2}-\mu_{\K}\right)} e^{n_{\K} t i\left(\log{2\pi} - \log{t} + 1 - \frac{\log{|\Delta_{\K}|}}{n_{\K}} + \frac{\log{w}}{n_{\K}}\right) + \beta(\mu_{\K},t)}}{z_1\cdots (z_1+n_{\K})} \left(1 + \alpha_{\K}(t)\right) \,dt
    \quad\text{and} \\
    J &= \int_{0}^{\infty} \frac{t^{n_{\K}\left(\frac{1}{2}-\nu_{\K}\right)} e^{n_{\K} t i\left(\log{2\pi} - \log{t} + 1 - \frac{\log{|\Delta_{\K}|}}{n_{\K}} + \frac{\log{w}}{n_{\K}}\right) + \beta(\nu_{\K},t)}}{z_2} \left(1 + \alpha_{\K}(t)\right) \,dt.
\end{align*}
These bounds are presented in the following result, whose proof is deferred until Section \ref{sec:ProofThmIJ}.

\begin{theorem}\label{thm:IJ}
If $g_{\K}^{\pm} = \frac{5}{8} + \frac{\pi}{2} \pm \frac{1}{n_{\K}} + \frac{3}{8 n_{\K}^2}$, then 
\begin{align*}
    |I| \leq 2^{\frac{n_{\K}}{2} + 2} e^{\frac{7}{2} n_{\K} + \frac{1}{442 n_{\K}} - \frac{1}{2}} \sqrt{g_{\K}^{+} n_{\K}}
    \qquad\text{and}\qquad
    |J| \leq 2^{\frac{n_{\K}}{2} + 2} e^{n_{\K}\left(\frac{1}{2} + \frac{4 n_{\K}}{n_{\K}-1} + \frac{1}{500 n_{\K}^6}\right) - \frac{1}{2}} \sqrt{g_{\K}^{-} n_{\K}}.
\end{align*}
\end{theorem}

\subsection{Proof of Theorem \ref{thm:Gertrude}}\label{sec:ProofGertrude}

Use \eqref{eqn:Kw_defn} and \eqref{eqn:Knw_defn} to obtain
\begin{align*}
    K(w)
    &= \left(\frac{(2\pi)^{\mu_{\K}}}{2}\right)^{n_{\K}} |\Delta_{\K}|^{\frac{1}{2}-\mu_{\K}} \frac{w^{\mu_{\K} + n_{\K}}}{2\pi i} \int_{-\infty}^{\infty} \frac{e^{n_{\K} t i\left(\log{2\pi} - \frac{\log{|\Delta_{\K}|}}{n_{\K}} + \frac{\log{w}}{n_{\K}}\right)}z_1^{-1}\cdots (z_1+n_{\K})^{-1}}{ \left(\cos{\frac{\pi z_1}{2}}\right)^{r_1 + r_2} \left(\sin{\frac{\pi z_1}{2}}\right)^{r_2} \Gamma(z_1)^{n_{\K}}} \,dt , \\
    K^{(n_{\K})}(w) 
    &= \left(\frac{(2\pi)^{\nu_{\K}}}{2}\right)^{n_{\K}} |\Delta_{\K}|^{\frac{1}{2}-\nu_{\K}} \frac{w^{\nu_{\K}}}{2\pi i} \int_{-\infty}^{\infty} \frac{e^{n_{\K} t i\left(\log{2\pi} - \frac{\log{|\Delta_{\K}|}}{n_{\K}} + \frac{\log{w}}{n_{\K}}\right)}}{z_2 \left(\cos{\frac{\pi z_2}{2}}\right)^{r_1 + r_2} \left(\sin{\frac{\pi z_2}{2}}\right)^{r_2} \Gamma(z_2)^{n_{\K}}} \,dt .
\end{align*}
Next, use \eqref{eqn:defsapplied} to see
\begin{align*}
    |K(w)|
    &\leq 2 (2\pi)^{\mu_{\K} n_{\K} - 1} |\Delta_{\K}|^{\frac{1}{2}-\mu_{\K}} w^{\mu_{\K} + n_{\K}} \!\left|\int_{0}^{\infty} \frac{e^{n_{\K} t i\left(\log{2\pi} - \frac{\log{|\Delta_{\K}|}}{n_{\K}} + \frac{\log{w}}{n_{\K}}\right) - \frac{n_{\K}\pi}{2}t}}{z_1\cdots (z_1+n_{\K}) \Gamma(z_1)^{n_{\K}}} \left(1 + \alpha_{\K}(t)\right) \,dt \right|, \\
    |K^{(n_{\K})}(w)|
    &\leq 2 (2\pi)^{\nu_{\K} n_{\K} - 1} |\Delta_{\K}|^{\frac{1}{2}-\nu_{\K}} w^{\nu_{\K}} \!\left|\int_{0}^{\infty} \frac{e^{n_{\K} t i\left(\log{2\pi} - \frac{\log{|\Delta_{\K}|}}{n_{\K}} + \frac{\log{w}}{n_{\K}}\right) - \frac{n_{\K}\pi}{2}t}}{z_2 \Gamma(z_2)^{n_{\K}}} \left(1 + \alpha_{\K}(t)\right) \,dt \right| .
\end{align*}
Further, use \eqref{eqn:GammaInsert} to see
\begin{align}
    |K(w)| &\leq 2 (2\pi)^{\mu_{\K} n_{\K} - 1} e^{-n_{\K} c} |\Delta_{\K}|^{\frac{1}{2}-\mu_{\K}} w^{\mu_{\K} + n_{\K}} |I| , \label{eqn:KwUpper}\\
    |K^{(n_{\K})}(w)| &\leq 2 (2\pi)^{\nu_{\K} n_{\K} - 1} e^{-n_{\K} c} |\Delta_{\K}|^{\frac{1}{2}-\nu_{\K}} w^{\nu_{\K}} |J|, \label{eqn:KnwUpper}
\end{align}
where $I$ and $J$ have been defined above Theorem \ref{thm:IJ}. Insert $|c|\leq 1/2$ and the bounds from Theorem \ref{thm:IJ} into \eqref{eqn:KwUpper} and \eqref{eqn:KnwUpper} to obtain Theorem \ref{thm:Gertrude}.

\section{Proof of Theorem \ref{thm:IJ}}\label{sec:ProofThmIJ}

Split $I$ and $J$ at $\mathfrak{s}_I$ and $\mathfrak{s}_J$ respectively, so that $I = I_1 + I_2$ and $J = J_1 + J_2$ such that
\begin{align*}
    I_1 &= \int_{0}^{\mathfrak{s}_I} \frac{t^{n_{\K}\left(\frac{1}{2}-\mu_{\K}\right)} e^{n_{\K} t i\left(\log{2\pi} - \log{t} + 1 - \frac{\log{|\Delta_{\K}|}}{n_{\K}} + \frac{\log{w}}{n_{\K}}\right) + \beta(\mu_{\K},t)}}{z_1\cdots (z_1+n_{\K})} \left(1 + \alpha_{\K}(t)\right) \,dt,\\
    I_2 &= \int_{\mathfrak{s}_I}^\infty \frac{t^{n_{\K}\left(\frac{1}{2}-\mu_{\K}\right)} e^{n_{\K} t i\left(\log{2\pi} - \log{t} + 1 - \frac{\log{|\Delta_{\K}|}}{n_{\K}} + \frac{\log{w}}{n_{\K}}\right) + \beta(\mu_{\K},t)}}{z_1\cdots (z_1+n_{\K})} \left(1 + \alpha_{\K}(t)\right) \,dt,
\end{align*}
and $J_1$, $J_2$ are defined similar to $I_1$, $I_2$ with $\mathfrak{s}_J$ in place of $\mathfrak{s}_I$, $\nu_{\K}$ in place of $\mu_{\K}$, and $z_2$ in place of $z_1$. The only conditions we require $\mathfrak{s}_I$ and $\mathfrak{s}_J$ to satisfy are that $|\beta(\sigma, t)| < 1$ for $t \geq \mathfrak{s}_I$ (when $\sigma = \mu_{\K}$) for $t \geq \mathfrak{s}_J$ (when $\sigma = \nu_{\K}$). These conditions ensure that we can use the Taylor series definition of $e^{\beta(\sigma,t)}$ to bound $I_2$ and $J_2$.

\subsection{Preliminary observations}\label{ssec:estimatefs_inv}

First, import \cite[Lem.~3.2.2]{SunleyThesis}, which establishes the bounds
\begin{align}
    \left|\beta(\sigma,t)\right| &\leq \frac{n_{\K}}{2 t} \left(3\sigma^2  + |\sigma| + \pi\right)
    \quad\text{for}\quad t\geq 0, \label{eqn:BetaInsert} \\
    \left|e^{\beta(\sigma,t)}\right| &\leq |s|^{n_{\K}\left(\frac{1}{2}-\sigma\right)} |t|^{n_{\K}\left(\sigma-\frac{1}{2}\right)} e^{n_{\K}\left(\sigma + \frac{2}{|\sigma|}\right)}.\label{eqn:ExpBetaInsert}
\end{align}
It follows that $|\beta(\mu_{\K},\mathfrak{s}_I)| < 1$ when $\mathfrak{s}_I > g_{\K}^{+} n_{\K}$ and $|\beta(\nu_{\K},\mathfrak{s}_J)| < 1$ when $\mathfrak{s}_J > g_{\K}^{-} n_{\K}$. We also need estimates for $\alpha_{\K}(t)$, which was defined underneath \eqref{eqn:defsapplied}.

\begin{lemma}\label{lem:fsinv_1}
If $\sigma\in\left\{\mu_{\K}, \nu_{\K}\right\}$, then the following statements are true:
\begin{itemize}
    \item If $t\geq 0$, then $|1 + \alpha_{\K}(t)| \leq 2^{\frac{n_{\K}}{2}}$.
    \item If $t>1$, then $|\alpha_{\K}(t)| \leq 1.09238\cdot 2^{n_{\K}}  e^{-\pi t}$.
\end{itemize}
\end{lemma}

\begin{proof}
To prove the result for $t\geq 0$, rewrite its definition in polar coordinate form. 
To prove the second statement, we refine the argument in \cite[Lem.~3.2.1]{SunleyThesis}. We will need the observations
\begin{equation*}
    \frac{1}{1+z_s} = 1 - z_s + \frac{{z_s}^2}{1+z_s}
    \quad\text{and}\quad
    \frac{1}{1-z_s} = 1 + z_s + \frac{{z_s}^2}{1-z_s}.
\end{equation*}
Now, note that $|1\pm z_s|$ satisfies
\begin{align*}
    \sqrt{e^{-2\pi t} \sin^2{\pi\sigma} + (1 \pm e^{-\pi t}\cos{\pi\sigma})^2} 
    = \sqrt{1 + e^{-2\pi t} \pm 2 e^{-\pi t}\cos{\pi\sigma}}
    &\geq \sqrt{1 + e^{-2\pi t} - 2 e^{-\pi t}}\\
    &= \sqrt{(1 - e^{-\pi t})^2} \\
    &= 1 - e^{-\pi t}.
\end{align*}
Therefore, $|1\pm z_s| \geq 1 - e^{-\pi}$ for $t\geq 1$. Next, use \eqref{eqn:NIST123} and the preceding to show that for $t\geq 1$,
\begin{align*}
    (1 + z_s)^{-\ell}
    = \left(1 - z_s + \frac{{z_s}^2}{1+z_s}\right)^{\ell}
    &= 1 + \sum_{k=1}^{\ell} (-1)^k \binom{\ell}{k} \left(z_s - \frac{{z_s}^2}{1+z_s}\right)^{k} \\
    &= 1 + O^*\!\!\left(\max_{1\leq k\leq\ell} \left|z_s - \frac{{z_s}^2}{1+z_s}\right|^{k} \left(\sum_{k=0}^{\ell} \binom{\ell}{k} - 1\right)\right) \\
    &= 1 + O^*\!\!\left((2^{\ell}-1) \max_{1\leq k\leq\ell} \left|\frac{z_s}{1+z_s}\right|^{k}\right) \\
    &= 1 + O^*\!\!\left((2^{\ell}-1) \max_{1\leq k\leq\ell} \left|\frac{z_s}{1-e^{-\pi}}\right|^{k}\right) \\
    &= 1 + O^*\!\!\left(\frac{(2^{\ell} - 1) |z_s|}{1 - e^{-\pi}}\right) .
\end{align*}
Following similar steps,
\begin{align*}
    (1 - z_s)^{-\ell}
    = \left(1 + z_s + \frac{{z_s}^2}{1-z_s}\right)^{\ell}
    &= 1 + \sum_{k=1}^{\ell} \binom{\ell}{k} \left(z_s + \frac{{z_s}^2}{1-z_s}\right)^{k} \\
    &= 1 + O^*\!\!\left(\max_{1\leq k\leq\ell} \left|\frac{z_s}{1-z_s}\right|^{k} \left(\sum_{k=0}^{\ell} \binom{\ell}{k} - 1\right)\right) \\
    &= 1 + O^*\!\!\left(\frac{(2^{\ell} - 1) |z_s|}{1 - e^{-\pi}}\right).
\end{align*}
It follows that
\begin{align*}
    \alpha_{\K}(t)
    &= \left(1 + O^*\!\!\left(\frac{(2^{r_2} - 1) e^{-\pi t}}{1 - e^{-\pi}}\right)\right) \left(1 + O^*\!\!\left(\frac{(2^{r_1+r_2} - 1) e^{-\pi t}}{1 - e^{-\pi}}\right)\right) - 1\\
    &= O^*\!\!\left(\frac{e^{-\pi t}}{1-e^{-\pi}}\left(2^{r_2} + 2^{r_1 + r_2} - 2 + \frac{(2^{n_{\K}} - 2^{r_2} - 2^{r_1 + r_2} + 1)e^{-\pi t}}{1-e^{-\pi}}\right)\right)\\
    &= O^*\!\!\left(\frac{e^{-\pi t}}{1-e^{-\pi}}\left(2^{r_2} + 2^{r_1 + r_2} - 2 + (2^{n_{\K}}+1) \frac{e^{-\pi}}{1-e^{-\pi}}\right)\right)\\
    &= O^*\!\!\left(\frac{e^{-\pi t}}{1-e^{-\pi}}\left(2^{n_{\K}}\left(1 + \frac{e^{-\pi}}{1-e^{-\pi}}\right) - 2 + \frac{e^{-\pi}}{1-e^{-\pi}} \right)\right)\\
    &= O^*\!\!\left(\frac{2^{n_{\K}} \left(1 + \frac{e^{-\pi}}{1-e^{-\pi}}\right)}{1-e^{-\pi}} e^{-\pi t}\right),
\end{align*}
using $\frac{e^{-\pi t}}{1-e^{-\pi}} < \frac{e^{-\pi}}{1-e^{-\pi}}$ for $t>1$, $-2 + \frac{e^{-\pi}}{1-e^{-\pi}} < 0$, and $2^{r_2} + 2^{r_1 + r_2} \leq 2^{n_{\K}}$.
\end{proof}

Next, we need the observation
\begin{equation}\label{eqn:Upsilon1PartialSummation}
    |\Upsilon_1(\nu_{\K}, u)| 
    \leq \frac{\xi_1(n_{\K})}{u},
\end{equation}
where
\begin{equation*}
    \Upsilon_1(\nu_{\K}, u) = \frac{ui}{n_{\K} \nu_{\K} + ui} - 1 = - \frac{n_{\K} \nu_{\K}}{n_{\K} \nu_{\K} + ui}
    \quad\text{and}\quad
    \xi_1(n_{\K}) = \frac{n_{\K} - 1}{2} .
\end{equation*}
Moreover, if $t$ is sufficiently large such that $|\beta(\sigma,t)| < 1$, then
\begin{equation}\label{eqn:ExpBetaTaylorSeries}
    \left|\sum_{k=1}^\infty \frac{\beta(\sigma,t)^{k}}{k!}\right|
    \leq \left|\frac{\beta(\sigma,t)}{1 - \beta(\sigma,t)}\right|
    := \beta_1(\sigma,t),
\end{equation}
and the following lemma provides estimates for $\beta_1(\sigma,t)$. 

\begin{lemma}\label{lem:BetaPartialSummation}
We have
\begin{align*}
    \beta_1\!\left(\mu_{\K},\frac{u}{n_{\K}}\right) &\leq \xi_2(n_{\K},\ell) \frac{n_{\K}^2}{u},
    \quad\text{where}\quad
    \xi_2(n_{\K},\ell) = g_{\K}^{+} \left|1 - \frac{g_{\K}^{+} n_{\K}^{2}}{\ell}\right|^{-1},\,\, u \geq \ell > g_{\K}^{+} n_{\K}^2 ,\\
    \beta_1\!\left(\nu_{\K},\frac{u}{n_{\K}}\right) &\leq \xi_3(n_{\K},\ell) \frac{n_{\K}^2}{u},
    \quad\text{where}\quad
    \xi_3(n_{\K},\ell) = g_{\K}^{-} \left|1 - \frac{g_{\K}^{-} n_{\K}^2}{\ell}\right|^{-1},\,\, u \geq \ell > g_{\K}^{-} n_{\K}^2 .
\end{align*}
\end{lemma}

\begin{proof}
Use \eqref{eqn:BetaInsert} to see that
\begin{align*}
    \left|\beta\!\left(\mu_{\K},\frac{u}{n_{\K}}\right)\right| 
    \leq \frac{n_{\K}^2 \left(3\mu_{\K}^2  + |\mu_{\K}| + \pi\right)}{2 u}
    = \frac{g_{\K}^{+} n_{\K}^2}{u}
    \quad\text{and}\quad
    \left|\frac{1}{1 - \beta\!\left(\mu_{\K},\frac{u}{n_{\K}}\right)}\right| 
    \leq \left|\frac{1}{1 - \frac{g_{\K}^{+} n_{\K}^2}{u}}\right| .
\end{align*}
It follows that
\begin{equation*}
    \beta_1\!\left(\mu_{\K},\frac{u}{n_{\K}}\right) 
    = \left|\frac{\beta\!\left(\mu_{\K},\frac{u}{n_{\K}}\right)}{1 - \beta\!\left(\mu_{\K},\frac{u}{n_{\K}}\right)}\right|
    < \frac{g_{\K}^{+} n_{\K}^2}{u} \left|1 - \frac{g_{\K}^{+} n_{\K}^2}{\ell}\right|^{-1}
    < \frac{\xi_2(n_{\K},\ell)  n_{\K}^2}{u} .
\end{equation*}
Follow a similar process with $\frac{1}{2} \left(3\nu_{\K}^2  + |\nu_{\K}| + \pi\right) = g_{\K}^{-}$ to establish the bound for $\beta_1(\nu_{\K},u/n_{\K})$.
\end{proof}

Another result we will need follows.

\begin{lemma}\label{lem:UpsilonPartialSummation}
If $u > 0$ and $\Upsilon(\mu_{\K}, u) = \prod_{j=0}^{n_{\K}} \frac{ui}{n_{\K} (\mu_{\K} + j) + ui} - 1$, 
then
\begin{equation*}
    |\Upsilon(\mu_{\K}, u)| \leq \left(\frac{n_{\K}+1}{2}\right)^{n_{\K}+1} \sum_{k=0}^{n_{\K}}  \binom{n_{\K}+1}{k} |n_{\K} \mu_{\K}|^{-k} u^{k-n_{\K}-1}.
\end{equation*}
\end{lemma}

\begin{proof}
Recall that $\mu_{\K} = -\tfrac{1}{2} - \tfrac{1}{2 n_{\K}}$. We have
\begin{align*}
    |\Upsilon(\mu_{\K}, u)|
    = \left|\frac{(u i)^{n_{\K}+1} - \prod_{j=0}^{n_{\K}} (n_{\K} (\mu_{\K} + j) + ui)}{\prod_{j=0}^{n_{\K}} (n_{\K} (\mu_{\K} + j) + ui)}\right|
    \leq \left|\frac{(u i)^{n_{\K}+1} - (n_{\K} \mu_{\K} + ui)^{n_{\K}+1}}{(ui)^{n_{\K}+1}}\right| .
\end{align*}
Using the binomial theorem and the triangle inequality, it follows that
\begin{align}
    |\Upsilon(\mu_{\K}, u)|
    &\leq \left|- (u i)^{-n_{\K}-1} \sum_{k=0}^{n_{\K}} \binom{n_{\K}+1}{k}(ui)^{k} (n_{\K} \mu_{\K})^{n_{\K}+1-k}\right| \nonumber\\
    &\leq u^{-n_{\K}-1} \sum_{k=0}^{n_{\K}} \left|\binom{n_{\K}+1}{k}(ui)^{k} (n_{\K} \mu_{\K})^{n_{\K}+1-k}\right| \nonumber\\
    &\leq \sum_{k=0}^{n_{\K}}  \binom{n_{\K}+1}{k} |n_{\K} \mu_{\K}|^{n_{\K}+1-k} u^{k-n_{\K}-1} . \qedhere
\end{align}
\end{proof}

Finally, we will need to estimate the integrals in the following lemma, using the theory of the incomplete $\Gamma$-function and the error function $\erf{\cdot}$. %
In particular, \cite[(8.4.6)]{NIST} tells us that
\begin{equation*}
    \Gamma(s,z) = \int_{z}^\infty y^{s-1} e^{-y} \,dy
    \quad\text{ satisfies}\quad 
    \Gamma\left(\frac{1}{2},z^2\right) = \pi^{\frac{1}{2}} \erfc{z},
\end{equation*}
in which $\Re{s} > 0$ and $\erfc{z}$ is the complement of the error function, defined as
\begin{equation*}
    \erfc{z} = \frac{2}{\sqrt{\pi}} \int_z^\infty e^{-t^2}\,dt ;
\end{equation*}
see \cite[\S 7.2(i)]{NIST} for more information on the error function.

\begin{lemma}\label{lem:IncompleteGammaLemma}
If $\ell >0$ and $k \in \{0,1,2\}$, then
\begin{align*}
    \int_{\ell}^\infty u^{-k-\frac{1}{2}} e^{-\frac{\pi u}{n_{\K}}} \,du
    \leq n_{\K}^{\frac{1}{2}} \ell^{-k} \erfc{n_{\K}^{-\frac{1}{2}} \pi^{\frac{1}{2}} \ell^{\frac{1}{2}}}.
\end{align*}
\end{lemma}

\begin{proof}
The result follows from $u^{-k} \leq \ell^{-k}$ and
\begin{equation*}
    \int_{\ell}^\infty u^{-\frac{1}{2}} e^{-\frac{\pi u}{n_{\K}}} \,du
    \leq \left(\frac{n_{\K}}{\pi}\right)^{\frac{1}{2}} \Gamma\!\left(\frac{1}{2}, \frac{\pi \ell}{n_{\K}} \right) 
    = n_{\K}^{\frac{1}{2}} \erfc{n_{\K}^{-\frac{1}{2}} \pi^{\frac{1}{2}} \ell^{\frac{1}{2}}}. \qedhere
\end{equation*}
\end{proof}

\subsection{Estimates for $I_1$ and $J_1$}

Use Lemma \ref{lem:fsinv_1}, $2/|\mu_{\K}| \leq 4$, and \eqref{eqn:ExpBetaInsert} to see that
\begin{align}
    |I_1| \leq 2^{\frac{n_{\K}}{2}} \int_{0}^{\mathfrak{s}_I} \frac{|z_1|^{n_{\K}\left(\frac{1}{2}-\mu_{\K}\right)} e^{n_{\K}\left(\mu_{\K} + \frac{2}{|\mu_{\K}|}\right)}}{|z_1|\cdots |z_1+n_{\K}|} \,dt
    &\leq 2^{\frac{n_{\K}}{2}} e^{\frac{7}{2} n_{\K} - \frac{1}{2}} \int_{0}^{\mathfrak{s}_I} \frac{|z_1|^{n_{\K} + \frac{1}{2}}}{|z_1|\cdots |z_1+n_{\K}|} \,dt \nonumber\\
    &\leq 2^{\frac{n_{\K}}{2}} e^{\frac{7}{2} n_{\K} - \frac{1}{2}} \int_{0}^{\mathfrak{s}_I} t^{- \frac{1}{2}} \,dt \nonumber\\
    &= 2^{\frac{n_{\K}}{2} + 1} e^{\frac{7}{2} n_{\K} - \frac{1}{2}} \sqrt{\mathfrak{s}_I} := \mathcal{I}_1(n_{\K}, \mathfrak{s}_I). \label{eqn:I1Estimate}
\end{align}
It follows from similar arguments, making changes \textit{mutatis mutandis}, that
\begin{equation}\label{eqn:J1Estimate}
    |J_1| \leq 2^{\frac{n_{\K}}{2} + 1} e^{n_{\K}\left(\frac{1}{2} + \frac{4 n_{\K}}{n_{\K}-1}\right) - \frac{1}{2}} \sqrt{\mathfrak{s}_J} := \mathcal{J}_1(n_{\K}, \mathfrak{s}_J) . 
\end{equation}

\subsection{Estimates for $I_2$ and $J_2$}\label{ssec:I2J2}

Use \eqref{eqn:ExpBetaTaylorSeries}, recall $w>0$, and apply the substitutions
\begin{equation*}
    u = n_{\K} t, \quad 1 + 2\log{v} = \log{n_{\K}} + 1 + \log{2\pi} + \frac{\log{w} - \log{|\Delta_{\K}|}}{n_{\K}},
\end{equation*}
to obtain
\begin{align}
    |I_2|
    &= n_{\K}^{-\frac{1}{2}} \left|\int_{\mathfrak{s}_I n_{\K}}^{\infty} \frac{u^{n_{\K}+\frac{1}{2}} e^{- u i\left(\log{u} -1 - 2\log{v} \right)}}{\prod_{j=0}^{n_{\K}} (n_{\K} (\mu_{\K} + j) + ui)} \sum_{\ell=0}^\infty \frac{\beta\!\left(\mu_{\K},\frac{u}{n_{\K}}\right)^\ell}{\ell!} \left(1 + \alpha_{\K}\!\left(\frac{u}{n_{\K}}\right)\right) \,du\right| \nonumber\\
    &\leq n_{\K}^{-\frac{1}{2}} \left|\int_{\mathfrak{s}_I n_{\K}}^{\infty} \frac{u^{-\frac{1}{2}} e^{- u i\left(\log{u} -1 - 2\log{v} \right)}}{\prod_{j=0}^{n_{\K}} \left(\frac{n_{\K} (\mu_{\K} + j) + u i}{u i}\right)} \sum_{\ell=0}^\infty \frac{\beta\!\left(\mu_{\K},\frac{u}{n_{\K}}\right)^{\ell}}{\ell!} \left(1 + \alpha_{\K}\!\left(\frac{u}{n_{\K}}\right)\right) \,du\right| \nonumber\\
    &\leq {n_{\K}}^{-\tfrac{1}{2}} \sum_{k=3}^{10} |I_k|, \nonumber 
\end{align}
in which
\begin{align*}
    I_3 &= \int_{\mathfrak{s}_I n_{\K}}^{\infty} u^{-\frac{1}{2}}e^{-u i\left(\log{u} -1 - 2\log{v} \right)} \,du,
    &&I_4 = \int_{\mathfrak{s}_I n_{\K}}^{\infty} u^{-\frac{1}{2}} \beta_1\!\left(\mu_{\K},\frac{u}{n_{\K}}\right) \,du, \\
    I_5 &= \int_{\mathfrak{s}_I n_{\K}}^{\infty} u^{-\frac{1}{2}} \alpha_{\K}\!\left(\frac{u}{n_{\K}}\right) \,du, 
    &&I_6 = \int_{\mathfrak{s}_I n_{\K}}^{\infty} u^{-\frac{1}{2}} \Upsilon(\mu_{\K}, u) \,du, \\
    I_7 &= \int_{\mathfrak{s}_I n_{\K}}^{\infty} u^{-\frac{1}{2}} \alpha_{\K}\!\left(\frac{u}{n_{\K}}\right) \beta_1\!\left(\mu_{\K},\frac{u}{n_{\K}}\right) \,du, 
    &&I_8 = \int_{\mathfrak{s}_I n_{\K}}^{\infty} u^{-\frac{1}{2}} \beta_1\!\left(\mu_{\K},\frac{u}{n_{\K}}\right) \Upsilon(\mu_{\K},u) \,du,\\
    I_9 &= \int_{\mathfrak{s}_I n_{\K}}^{\infty} u^{-\frac{1}{2}} \alpha_{\K}\!\left(\frac{u}{n_{\K}}\right) \Upsilon(\mu_{\K},u) \,du, 
    &&I_{10} = \int_{\mathfrak{s}_I n_{\K}}^{\infty} u^{-\frac{1}{2}} \alpha_{\K}\!\left(\frac{u}{n_{\K}}\right) \Upsilon(\mu_{\K},u) \beta_1\!\left(\mu_{\K},\frac{u}{n_{\K}}\right) du .
\end{align*}
Define $J_i$ similarly to $I_i$ for $3\leq i\leq 10$, with $\mathfrak{s}_J$ in place of $\mathfrak{s}_I$, $\nu_{\K}$ in place of $\mu_{\K}$, and $\Upsilon_1$ in place of $\Upsilon$. Using similar arguments, we see that
\begin{equation*}
    |J_2| \leq n_{\K}^{-\frac{1}{2}} \sum_{k=3}^{10} |J_k|.
\end{equation*}

To estimate $I_3$ and $J_3$, we use \cite[Satz 204]{LandauEinfuehrung}, which states that for $U>0$ and $\mu\in\mathbb{R}$, 
\begin{equation*}
    \left| \int_{0}^{U} u^{-\frac{1}{2}} e^{-ui (\log{u} - \mu)} \,du \right| < 26.
\end{equation*}
Let $U\to\infty$ and note that $1 + 2\log{v}$ is real, so it follows that
\begin{equation}\label{eqn:I3J3estimate}
    |I_3|, |J_3| < 26.
\end{equation}
Since $w>0$ and $|\Delta_{\K}|$ can be large, there are no clear bounds one could place on $1 + 2\log{v}$, nor would there be significant improvements (if any) available from considering $\mu$ in a restricted range.

To estimate $I_i$ for $4\leq i\leq 10$, use Lemmas \ref{lem:fsinv_1}-\ref{lem:IncompleteGammaLemma} to obtain
\begin{align}
    |I_4| &\leq \xi_2(n_{\K},\mathfrak{s}_I n_{\K}) n_{\K}^2 \int_{\mathfrak{s}_I n_{\K}}^\infty u^{-\frac{3}{2}}\,du
    = 2 n_{\K}^2 \frac{\xi_2(n_{\K},\mathfrak{s}_I n_{\K})}{\sqrt{\mathfrak{s}_I n_{\K}}} , \label{eqn:I4Estimate}\\
    |I_5| &\leq 1.1 n_{\K}^{\frac{1}{2}} 2^{n_{\K}} \erfc{\sqrt{\pi \mathfrak{s}_I}} , \label{eqn:I5Estimate}\\
    |I_6| &\leq 
    \left(\frac{n_{\K}+1}{2}\right)^{n_{\K}+1} \sum_{k=0}^{n_{\K}}  \binom{n_{\K}+1}{k} |n_{\K} \mu_{\K}|^{-k} \int_{\mathfrak{s}_I n_{\K}}^\infty u^{k-n_{\K}-\frac{3}{2}}\,du \nonumber\\
    &=2 \left(\frac{n_{\K}+1}{2}\right)^{n_{\K}+1} \sum_{k=0}^{n_{\K}}  \binom{n_{\K}+1}{k} \frac{|n_{\K} \mu_{\K}|^{-k} (\mathfrak{s}_I n_{\K})^{k-n_{\K}-\frac{1}{2}}}{2(n_{\K}-k)+1} , \label{eqn:I6Estimate}\\
    |I_7| &\leq \frac{1.1}{\mathfrak{s}_I} \xi_2(n_{\K},\mathfrak{s}_I n_{\K}) n_{\K}^{\frac{3}{2}} 2^{n_{\K}} \erfc{\sqrt{\pi \mathfrak{s}_I}}, \label{eqn:I7Estimate}\\
    |I_8| &\leq 
    \xi_2(n_{\K}, \mathfrak{s}_I n_{\K}) n_{\K}^2 \left(\frac{n_{\K}+1}{2}\right)^{n_{\K}+1} \sum_{k=0}^{n_{\K}}  \binom{n_{\K}+1}{k} |n_{\K} \mu_{\K}|^{-k} \int_{\mathfrak{s}_I n_{\K}}^\infty u^{k-n_{\K}-\frac{5}{2}}\,du \nonumber\\
    &= 2 \xi_2(n_{\K}, \mathfrak{s}_I n_{\K}) n_{\K}^2 \left(\frac{n_{\K}+1}{2}\right)^{n_{\K}+1} \sum_{k=0}^{n_{\K}}  \binom{n_{\K}+1}{k} \frac{|n_{\K} \mu_{\K}|^{-k} (\mathfrak{s}_I n_{\K})^{k-n_{\K}-\frac{3}{2}}}{3+2(n_{\K}-k)} , \label{eqn:I8Estimate}\\
    |I_9| &\leq \frac{1.1}{\mathfrak{s}_I} n_{\K}^{-\frac{1}{2}} 2^{n_{\K}} \xi_4(n_{\K}) \erfc{\sqrt{\pi \mathfrak{s}_I}} , \label{eqn:I9Estimate}\\
    |I_{10}| &\leq \frac{1.1}{\mathfrak{s}_I^2} n_{\K}^{\frac{1}{2}} 2^{n_{\K}} \xi_2(n_{\K},\mathfrak{s}_I n_{\K}) \xi_4(n_{\K}) \erfc{\sqrt{\pi \mathfrak{s}_I}}. \label{eqn:I10Estimate}
\end{align}
If $\mathfrak{s}_I$ is sufficiently large, then the presence of $\erfc{\cdot}$ in \eqref{eqn:I5Estimate}, \eqref{eqn:I7Estimate}, \eqref{eqn:I9Estimate}, and \eqref{eqn:I10Estimate} means that $I_5$, $I_7$, $I_9$, and $I_{10}$ contribute a negligible amount to the upper estimate for $|I_2|$, because $\erfc{w} \to 0$ very quickly for $w>0$. In fact, a simple upper bound $\erfc{x} \leq \exp(-x^2)$ is found in \cite{ChangCosmanMilstein}. It follows that an upper estimate for $|I_2|$ is mainly derived from \eqref{eqn:I4Estimate}, \eqref{eqn:I6Estimate}, and \eqref{eqn:I8Estimate}. Moreover, if $\mathfrak{s}_I = \tau_I g_{\K}^{+} n_{\K}$ such that $1< \tau_I < 1.25$, then the sum of the upper bounds in \eqref{eqn:I5Estimate}, \eqref{eqn:I7Estimate}, \eqref{eqn:I9Estimate}, \eqref{eqn:I10Estimate} decreases as $n_{\K}$ increases, so there is a constant $c_I$ such that
\begin{align*}
    \sum_{i\in\{5,7,9,10\}} |I_i| < c_I.
\end{align*}
In particular, if $\tau_I = e^{\frac{1}{221 n_{\K}}}$, then $c_I = 10^{-5}$ is admissible; this is obtained by computing the sum of the upper bounds in \eqref{eqn:I5Estimate}, \eqref{eqn:I7Estimate}, \eqref{eqn:I9Estimate}, \eqref{eqn:I10Estimate} at $n_{\K} = 2$. We choose this definition of $\mathfrak{s}_I$ to ensure the necessary condition $\mathfrak{s}_I > g_{\K}^{+} n_{\K}$ is satisfied. Now, combine observations to see
\begin{equation}\label{eqn:checkcondition}
    \sum_{k=4}^{10} |I_k| \leq 6 \max\left\{\xi_2(n_{\K},\mathfrak{s}_I n_{\K}) \frac{n_{\K}^{3/2}}{\sqrt{\mathfrak{s}_I}} , \xi_4(n_{\K}, \mathfrak{s}_I), \xi_2(n_{\K}, \mathfrak{s}_I n_{\K}) \xi_5(n_{\K}, \mathfrak{s}_I) n_{\K}^2 \right\} + c_I ,
\end{equation}
in which
\begin{align*}
    \xi_4(n_{\K}, \mathfrak{s}_I) &= \left(\frac{n_{\K}+1}{2}\right)^{n_{\K}+1} \sum_{k=0}^{n_{\K}}  \binom{n_{\K}+1}{k} \frac{|n_{\K} \mu_{\K}|^{-k} (\mathfrak{s}_I n_{\K})^{k-n_{\K}-\frac{1}{2}}}{2(n_{\K}-k)+1},\\
    \xi_5(n_{\K}, \mathfrak{s}_I) &= \left(\frac{n_{\K}+1}{2}\right)^{n_{\K}+1} \sum_{k=0}^{n_{\K}}  \binom{n_{\K}+1}{k} \frac{|n_{\K} \mu_{\K}|^{-k} (\mathfrak{s}_I n_{\K})^{k-n_{\K}-\frac{3}{2}}}{3+2(n_{\K}-k)} .
\end{align*}

Similarly, use Lemmas \ref{lem:fsinv_1}-\ref{lem:IncompleteGammaLemma} and \eqref{eqn:Upsilon1PartialSummation} to obtain
\begin{align}
    |J_4| &\leq 2 n_{\K}^2 \frac{\xi_3(n_{\K},\mathfrak{s}_J n_{\K})}{\sqrt{\mathfrak{s}_J n_{\K}}} , \nonumber\\ 
    |J_5| &\leq 1.1 n_{\K}^{\frac{1}{2}} 2^{n_{\K}} \erfc{\sqrt{\pi \mathfrak{s}_J}} , \nonumber\\ 
    |J_6| &\leq \frac{2\,\xi_1(n_{\K})}{\sqrt{\mathfrak{s}_J n_{\K}}} , \nonumber\\ 
    |J_7| &\leq \frac{1.1}{\mathfrak{s}_J} \xi_3(n_{\K},\mathfrak{s}_J n_{\K}) n_{\K}^{\frac{3}{2}} 2^{n_{\K}} \erfc{\sqrt{\pi \mathfrak{s}_J}}, \nonumber\\ 
    |J_8| &\leq \frac{2\,\xi_3(n_{\K}, \mathfrak{s}_J n_{\K}) \xi_1(n_{\K}) n_{\K}^2}{3 (\mathfrak{s}_J n_{\K})^{3/2}} , \nonumber\\ 
    |J_9| &\leq \frac{1.1}{\mathfrak{s}_J} n_{\K}^{-\frac{1}{2}} 2^{n_{\K}} \xi_1(n_{\K}) \erfc{\sqrt{\pi \mathfrak{s}_J}} , \nonumber\\ 
    |J_{10}| &\leq \frac{1.1}{\mathfrak{s}_J^2} n_{\K}^{\frac{1}{2}} 2^{n_{\K}} \xi_3(n_{\K},\mathfrak{s}_J n_{\K}) \xi_1(n_{\K}) \erfc{\sqrt{\pi \mathfrak{s}_J}}. \nonumber 
\end{align}
Moreover, if $\mathfrak{s}_J = \tau_J g_{\K}^{-} n_{\K}$ such that $1< \tau_J < 1.25$, then there is a constant $c_J$ such that
\begin{align*}
    \sum_{i\in\{5,7,9,10\}} |J_i| < c_J,
\end{align*}
and $\tau_J = e^{\frac{1}{250 n_{\K}^5}}$ implies that $c_J = 10^{-2}$ is admissible; this is obtained by summing the appropriate bounds at $n_{\K} = 2$. Combine these observations as before to see
\begin{equation}\label{eqn:checkcondition1}
    \sum_{k=4}^{10} |J_k| 
    \leq 6 \max\left\{ \frac{\xi_3(n_{\K},\mathfrak{s}_J n_{\K}) n_{\K}^{3/2}}{\sqrt{\mathfrak{s}_J}}, \frac{\xi_1(n_{\K})}{\sqrt{\mathfrak{s}_J n_{\K}}}, \frac{\xi_3(n_{\K}, \mathfrak{s}_J n_{\K}) \xi_1(n_{\K}) \sqrt{n_{\K}}}{3 \mathfrak{s}_J^{3/2}} \right\} + c_J.
\end{equation}

Finally, let $\mathfrak{s}_I = \tau_I g_{\K}^{+} n_{\K}$ and combine \eqref{eqn:I3J3estimate} with \eqref{eqn:checkcondition} to obtain
\begin{equation}\label{eqn:I2Estimate}
    |I_2| 
    \leq \underbrace{\frac{26}{\sqrt{n_{\K}}} + c_I + 6 \max\left\{ \frac{\xi_2(n_{\K},\mathfrak{s}_I n_{\K})}{\sqrt{\mathfrak{s}_I}} n_{\K}, \frac{\xi_4(n_{\K}, \mathfrak{s}_I)}{\sqrt{n_{\K}}}, \xi_2(n_{\K}, \mathfrak{s}_I n_{\K}) \xi_5(n_{\K}, \mathfrak{s}_I) n_{\K}^{3/2} \right\}}_{\mathcal{I}_2(n_{\K}, \mathfrak{s}_I)} .
\end{equation}
Similarly, let $\mathfrak{s}_J = \tau_J g_{\K}^{-} n_{\K}$ and combine \eqref{eqn:I3J3estimate} with \eqref{eqn:checkcondition1} to obtain
\begin{equation}\label{eqn:J2Estimate}
    |J_2| 
    \leq \underbrace{\frac{26}{\sqrt{n_{\K}}} + c_J + 6 \max\left\{ \frac{\xi_3(n_{\K},\mathfrak{s}_J n_{\K})}{\sqrt{\mathfrak{s}_J}} n_{\K}, \frac{\xi_1(n_{\K})}{\sqrt{\mathfrak{s}_J} n_{\K}}, \frac{\xi_3(n_{\K}, \mathfrak{s}_J n_{\K}) \xi_1(n_{\K})}{3 \mathfrak{s}_J^{3/2}} \right\}}_{\mathcal{J}_2(n_{\K}, \mathfrak{s}_J)}.
\end{equation}

\subsection{Proof of Theorem \ref{thm:IJ}}

\begin{table}[]
    \centering
    \begin{tabular}{ccccc}
        $n_{\K}$ & $\mathcal{I}_1(n_{\K}, \mathfrak{s}_I)$ & $\mathcal{I}_2(n_{\K}, \mathfrak{s}_I)$ & $\mathcal{I}_1(n_{\K}, \mathfrak{s}_I)+\mathcal{I}_2(n_{\K}, \mathfrak{s}_I)$ & $2\, \mathcal{I}_1(n_{\K}, \mathfrak{s}_I)$ \\
        \hline
        2  &  $6.29140\cdot 10^{3}$  &  $6.28243\cdot 10^{3}$  &  $1.25738\cdot 10^{4}$  &  $1.25828\cdot 10^{4}$  \\
        3  &  $3.46291\cdot 10^{5}$  &  $1.10624\cdot 10^{4}$  &  $3.57354\cdot 10^{5}$  &  $6.92583\cdot 10^{5}$  \\
        4  &  $1.83494\cdot 10^{7}$  &  $1.66822\cdot 10^{4}$  &  $1.83661\cdot 10^{7}$  &  $3.66988\cdot 10^{7}$  \\
        5  &  $9.49233\cdot 10^{8}$  &  $2.30302\cdot 10^{4}$  &  $9.49256\cdot 10^{8}$  &  $1.89847\cdot 10^{9}$  \\
        10  &  $2.95297\cdot 10^{17}$  &  $6.35947\cdot 10^{4}$  &  $2.95297\cdot 10^{17}$  &  $5.90594\cdot 10^{17}$  \\
        20  &  $2.09477\cdot 10^{34}$  &  $1.77778\cdot 10^{5}$  &  $2.09477\cdot 10^{34}$  &  $4.18954\cdot 10^{34}$  \\
        30  &  $1.29704\cdot 10^{51}$  &  $3.25342\cdot 10^{5}$  &  $1.29704\cdot 10^{51}$  &  $2.59409\cdot 10^{51}$  \\
    \end{tabular}
    \caption{Comparison between $\mathcal{I}_1(n_{\K}, \mathfrak{s}_I)$ and $\mathcal{I}_2(n_{\K}, \mathfrak{s}_I)$ for $\mathfrak{s}_I = g_{\K}^{+} e^{\frac{1}{221 n_{\K}}} n_{\K}$ and several $n_{\K}$.}
    \label{tab:IEstimateUpper}
\end{table}

\begin{table}[]
    \centering
    \begin{tabular}{ccccc}
        $n_{\K}$ & $\mathcal{J}_1(n_{\K}, \mathfrak{s}_J)$ & $\mathcal{J}_2(n_{\K}, \mathfrak{s}_J)$ & $\mathcal{J}_1(n_{\K}, \mathfrak{s}_J)+\mathcal{J}_2(n_{\K}, \mathfrak{s}_J)$ & $2\, \mathcal{J}_1(n_{\K}, \mathfrak{s}_J)$ \\
        \hline
        2  &  $1.10875\cdot 10^{8}$  &  $9.08271\cdot 10^{4}$  &  $1.10966\cdot 10^{8}$  &  $2.21749\cdot 10^{8}$  \\
        3  &  $2.41314\cdot 10^{9}$  &  $8.71192\cdot 10^{5}$  &  $2.41401\cdot 10^{9}$  &  $4.82628\cdot 10^{9}$  \\
        4  &  $1.85208\cdot 10^{11}$  &  $4.31093\cdot 10^{6}$  &  $1.85212\cdot 10^{11}$  &  $3.70416\cdot 10^{11}$  \\
        5  &  $1.90865\cdot 10^{13}$  &  $1.48631\cdot 10^{7}$  &  $1.90865\cdot 10^{13}$  &  $3.81729\cdot 10^{13}$  \\
        10  &  $5.29109\cdot 10^{23}$  &  $6.87312\cdot 10^{8}$  &  $5.29109\cdot 10^{23}$  &  $1.05822\cdot 10^{24}$  \\
        20  &  $6.69420\cdot 10^{44}$  &  $3.14518\cdot 10^{10}$  &  $6.69420\cdot 10^{44}$  &  $1.33884\cdot 10^{45}$  \\
        30  &  $8.55550\cdot 10^{65}$  &  $2.93612\cdot 10^{11}$  &  $8.55550\cdot 10^{65}$  &  $1.71110\cdot 10^{66}$  \\
    \end{tabular}
    \caption{Comparison between $\mathcal{J}_1(n_{\K}, \mathfrak{s}_J)$ and $\mathcal{J}_2(n_{\K}, \mathfrak{s}_J)$ for $\mathfrak{s}_J = e^{\frac{1}{250 n_{\K}^5}} g_{\K}^{-} n_{\K}$ and several $n_{\K}$.}
    \label{tab:JEstimateUpper}
\end{table}

Use \eqref{eqn:I1Estimate}, \eqref{eqn:I2Estimate} to see $|I| \leq |I_1| + |I_2| \leq \mathcal{I}_1(n_{\K}, \mathfrak{s}_I) + \mathcal{I}_2(n_{\K}, \mathfrak{s}_I)$ when $\mathfrak{s}_I = \tau_I g_{\K}^{+} n_{\K}$ and \eqref{eqn:J1Estimate}, \eqref{eqn:J2Estimate} to see $|J| \leq \mathcal{J}_1(n_{\K}, \mathfrak{s}_J) + \mathcal{J}_2(n_{\K}, \mathfrak{s}_J)$ when $\mathfrak{s}_J = \tau_J g_{\K}^{-} n_{\K}$. 
Now, if $\mathfrak{s}_I$ and $\mathfrak{s}_J$ are chosen appropriately, then we also have
\begin{align}
    \mathcal{I}_1(n_{\K}, \mathfrak{s}_I) + \mathcal{I}_2(n_{\K}, \mathfrak{s}_I) &\leq 2 \mathcal{I}_1(n_{\K}, \mathfrak{s}_I)
    \quad\text{and} \label{eqn:bbt1}\\
    \mathcal{J}_1(n_{\K}, \mathfrak{s}_J) + \mathcal{J}_2(n_{\K}, \mathfrak{s}_J) &\leq 2 \mathcal{J}_1(n_{\K}, \mathfrak{s}_J). \label{eqn:bbt2}
\end{align}
To this end, we will choose
\begin{equation*}
    \mathfrak{s}_I = e^{\frac{1}{221 n_{\K}}} g_{\K}^{+} n_{\K}
    \qquad\text{and}\qquad 
    \mathfrak{s}_J = e^{\frac{1}{250 n_{\K}^5}} g_{\K}^{-} n_{\K},
\end{equation*}
because smaller choices of $\mathfrak{s}_I$ and $\mathfrak{s}_J$ are favourable and these choices satisfy \eqref{eqn:bbt1}, \eqref{eqn:bbt2} for all $n_{\K} \geq 2$; evidence to this end is presented in Tables \ref{tab:IEstimateUpper} and \ref{tab:JEstimateUpper}. To obtain the result, insert these choices of $\mathfrak{s}_I$ and $\mathfrak{s}_J$ into \eqref{eqn:bbt1} and \eqref{eqn:bbt2}.


\begin{remark}
The choice for $\tau_J$ in $\mathfrak{s}_J$ was arbitrary, since any choice satisfying $\mathfrak{s}_J > g_{\K}^{-} n_{\K}$ also satisfies \eqref{eqn:bbt2}; the appeal of $\exp(1/(250 n_{\K}^5))$ is its quick convergence rate. The choice for $\tau_I$ in $\mathfrak{s}_I$ was \textit{not} arbitrary, since any choice $\mathfrak{s}_I > g_{\K}^{+} n_{\K}$ satisfies \eqref{eqn:bbt1} when $n_{\K} \geq 4$, but larger choices are needed for $n_{\K}\in\{2,3\}$. To determine the choice for $\mathfrak{s}_I$, we wrote $\mathfrak{s}_I = \exp(1/(c\, n_{\K})) g_{\K}^{+} n_{\K}$ and searched for the largest integer $c$ such that \eqref{eqn:bbt1} holds for $n_{\K} \geq 2$. 
\end{remark}

\bibliographystyle{amsplain}
\bibliography{references}

\end{document}